\documentclass[10pt,twoside]{siamart1116}

\usepackage[english]{babel}
\usepackage{graphicx,epstopdf,epsfig}
\usepackage{amsfonts,epsfig,fancyhdr,graphics, hyperref,amsmath,amssymb}

\usepackage{dsfont}
\usepackage[shortlabels]{enumitem}
\usepackage[textsize=footnotesize]{todonotes}
\usepackage{pgfplots}
\usepackage{tikz}

\def\cvd{~\vbox{\hrule\hbox{%
     \vrule height1.3ex\hskip0.8ex\vrule}\hrule } }

\newcommand{\Names}{Zachary Brennan}
\newcommand{\Title}{Probabilistic zero forcing with vertex reversion}

\newtheorem{remark}[theorem]{Remark}

\newtheorem{fact}[theorem]{Fact}
\newtheorem{observation}[theorem]{Observation}

\newcommand{\reals}{\mathbb{R}}

\def\E{{\mathbf{E}}}
\def\P{{\mathbf{P}}}
\def\varep{{\varepsilon}}
\newcommand{\naturals}{\mathbb{N}}
\newcommand{\defn}[1]{\emph{#1}}
\newcommand{\bo}[1]{\mathbf{#1}}
\newcommand{\mf}[1]{\mathfrak{#1}}
\DeclareMathOperator{\probabs}{pta}
\DeclareMathOperator{\expabs}{eta}
\newcommand{\norm}[1]{\left\Vert #1\right\Vert}
\renewcommand{\exp}[1]{\operatorname{exp}\left\{#1\right\}}
 
\renewcommand{\theequation}{\thesection.\arabic{equation}}
\renewtheorem{theorem}{Theorem}[section]

\begin{document}

\bibliographystyle{plain}

\setcounter{page}{1}

\thispagestyle{empty}

 \title{\Title}

\author{
Zachary Brennan\thanks{Department of Mathematics,
Iowa State University, Ames, Iowa 50011, USA
(brennanz@iastate.edu).}
}

\markboth{\Names}{\Title}

\maketitle

\begin{abstract}
Probabilistic zero forcing is a graph coloring process in which blue vertices ``infect" (color blue) white vertices with a probability proportional to the number of neighboring blue vertices. 
We introduce reversion probabilistic zero forcing (RPZF), which shares the same infection dynamics but also allows for blue vertices to revert to being white in each round. 
We establish a tool which, given a graph's RPZF Markov transition matrix, calculates the probability that the graph turns all white or all blue as well as the time at which this is expected to occur.
For specific graph families we produce a threshold number of blue vertices for the graph to become entirely blue in the next round with high probability.
\end{abstract}




\section{Introduction}
Zero forcing is a graph coloring process which was introduced independently as a condition for the control of quantum systems \cite{Quantum1} and as a bound for the maximum nullity of a matrix in the study of the minimum rank problem \cite{AIM_Session_2008}. 
Zero forcing has since been used extensively in the study of the minimum rank problem (see \cite{HandbookZF} and the references therein), and has been found to have further connections with graph search algorithms \cite{YANG13}, power domination \cite{BFFFHVW18}, and the Cops and Robbers game \cite{BREEN18}. 
These connections have led to the study of zero forcing in its own right, and variants of zero forcing have since emerged (see the workshop summary \cite{ZFVariants} for examples).
This paper focuses on probabilistic variations of zero forcing, the first of which was introduced by Kang and Yi in \cite{KY13}.
Specifically, in this paper we introduce reversion probabilistic zero forcing (RPZF), a process in which blue vertices can now revert back to being white. 
This process can be formulated as a Markov chain and is thus studied using the theory of Markov chains.
Probabilistic zero forcing cannot move across graph components and is typically studied on connected graphs. Thus, throughout this paper we assume $G$ is a simple connected graph on $n$ vertices.

RPZF can also be viewed as a discrete-time analog of the susceptible-infected-susceptible (SIS) contact process, a basic continuous-time model for the spread of infection. 
In particular, RPZF is an example of a semi-heterogeneous model, where the infection rate is different for different vertices but the recovery rate is uniform. 
Contact processes are traditionally considered on infinite lattices with the goal of understanding long-term behavior, though more recent work has also been done on other large graph structures such as scale-free and power-law graphs \cite{gomez2010discrete, jacob2017contact}.
In contrast, we initiate the study of RPZF by considering its long-term behavior on the complete and complete bipartite graphs, which model denser and more interconnected populations. 
Additionally, most discretized SIS models only consider distance-1 interactions, where the probability of a vertex being infected depends only on its neighborhood \cite{pare2018analysis}. 
By contrast, RPZF is a distance-2 model where an infected vertex is \textit{more infectious} when it has more infected neighbors. 

\Cref{sec:RPZF} defines RPZF then formulates it as a Markov chain, and
detailed discussion of how RPZF relates to the SIS contact process is provided in \cref{sec:contact}.
Probability theory, in particular the theory of Markov chains, is used in \cref{sec:gen_rpzf} to calculate RPZF parameters, characterize each state of the RPZF Markov chain as transient or absorbing, and make assertions about the exit probabilities of RPZF. 
These results enable the quantification of RPZF behavior on any finite graph, given its RPZF-matrix representation.
\Cref{sec:thresholds} goes into detailed analysis of RPZF on the complete graph, complete balanced bipartite graph, and star graph. In particular, threshold-like results are developed to prove the required number of blue vertices to force the entire graph blue in one step with high probability.
Finally, simulations and numerical approximations of RPZF and its parameters are provided for a number of graph families in \cref{sec:sims}.

\subsection{Basic Notation}
A \defn{graph} is a pair $G=(V,E)$ where the set $E=E(G)$ of \defn{edges} consists of 2-element subsets of $V=V(G)$, the finite set of \defn{vertices}. Thus all graphs discussed (except in \cref{sec:contact}) are simple, undirected, and finite. Two vertices $v,w\in V$ are \defn{adjacent} if $\{v,w\}\in E$. The \defn{open neighborhood of $v$} is the set of all vertices adjacent to $v$, denoted by
\(N(v)=\{w\in V:\{v,w\}\in E\}.\) The \defn{degree of $v$} is $\deg v=|N(v)|$, the number of vertices adjacent to $v$, and the \defn{closed neighborhood of $v$} is $N[v]=N(v)\cup\{v\}$. We say $G$ is \defn{connected} if for every $v,w\in V$ there exists a \defn{path} of vertices $v=v_0,v_1,\ldots,v_k=w$ such that $\{v_{i-1},v_i\}\in E$ for all $i\in\{1,\ldots,k\}$. In this paper we consider only connected graphs.

If $A=[a_{ij}]$ is a matrix, then $a_{ij}^{(k)}$ and $(A^k)_{ij}$ both denote the $(i,j)$th entry of $A^k$. The identity matrix is denoted by $I$, and the matrix containing all zero entries is denoted by $O$. We will use $\mathds{1}=[1~1~\cdots~1]^T$ to denote the column vector containing all ones and $\naturals=\{0,1,2,\ldots\}$ to refer to the set of non-negative integers.

\subsection{Zero Forcing}
Suppose $G$ is colored so that every vertex is blue or white. The \defn{(deterministic) zero forcing color change rule} describes how the vertices of $G$ can change color: a blue vertex $u$ will \defn{force} (change) a white vertex $w$ to be blue if $w$ is the only white neighbor of $u$. This is denoted by $u\rightarrow w$. 

Probabilistic zero forcing (PZF) is an extension of the deterministic model. Let $B\subseteq V(G)$ be a set of blue vertices. 
In one round of PZF, each blue vertex $u\in B$ attempts to force each of its white neighbors $w\in N(u)\setminus B$ independently with probability 
\[\P (u\to w)=\frac{|N[u]\cap B|}{\deg{u}}.\]

In deterministic zero forcing, a parameter of interest is the minimum number of vertices required to eventually force the entire graph blue. 
In probabilistic zero forcing, this parameter is trivial because only one vertex is required in a connected graph. 
Thus, the study of probabilistic zero forcing is instead primarily concerned with how long it takes the entire graph to be forced blue starting from one blue vertex \cite{GH18}.

The \defn{probabilistic propagation time of $B$}, denoted $\operatorname{pt}_{pzf}(G,B)$, is the first round in which all vertices are blue, starting from the initial set of blue vertices $B$ \cite{GH18}. Note that $\operatorname{pt}_{pzf}(G,B)$ is a random variable and a \defn{round} constitutes a single application of the probabilistic color change rule and that the first application occurs in round 1. The \defn{expected propagation time of $B$} is then defined as \cite{GH18}
\[\operatorname{ept}(G,B)=\E[\operatorname{pt}_{pzf}(G,B)].\]
The \defn{expected propagation time} of a connected graph $G$ is the minimum expected propagation time of $B$ over all one-vertex sets $B$ of $G$ \cite{GH18}, i.e.,
\[\operatorname{ept(G)}=\min\{\operatorname{ept}(G,\{v\}):v\in V(G)\}.\]

\newpage

\subsection{The Markov approach to probabilistic zero forcing}
Probabilistic zero forcing can also be formulated as a time-homogeneous Markov chain.
The study of PZF with Markov chains was introduced in \cite{KY13} and studied further in \cite{CC+20,GH18}.
Let $X_0,X_1,\ldots$ be random variables that take on values from the finite state space $\mathcal{S}=\{S_0,S_1,\ldots,S_s\}$. 
We say $(X_t)$ is a \defn{time-homogeneous (discrete-time) Markov chain} on the state space $\mathcal{S}$ if for any time $t\in\naturals$ and any states $S_j,S_i,S_{i_{t-1}},S_{i_{t-2}},\ldots,S_{i_0}\in\mathcal{S}$,
\begin{equation}\label{eqn:markovprop}
    \P[X_{t+1}=S_j~|~X_t=S_i,X_{t-1}=S_{i_{t-1}},\ldots,X_0=S_{i_0}]=\P[X_{t+1}=S_j~|~X_t=S_i]
\end{equation}
and
\begin{equation}\label{eqn:time-hom}
\P[X_{t+1}=S_j~|~X_t=S_i]=\P[X_1=S_j~|~X_0=S_i]
\end{equation}
\Cref{eqn:markovprop} is the \defn{Markov property} and it says that given the current state $X_t$, any other information about the history of the chain is irrelevant for predicting $X_{t+1}$. \Cref{eqn:time-hom} is the definition of time-homogeneity, and it says that the transition probability does not depend on the current time $t$.
We will consider only time-homogeneous Markov chains.

The \defn{Markov transition matrix} $M=[m_{ij}]$ is the $(s+1)\times (s+1)$ matrix defined by $$m_{ij}=\P[X_1=S_j~|~X_0=S_i],$$ the probability of the chain going from state $S_i$ to state $S_j$ in one time-step. Markov transition matrices are \textit{right stochastic}, i.e., $\sum_{k=0}^{s}m_{ik}=1$ for all $i=0,1,\ldots,s$. 
It follows immediately from the Markov property and time-homogeneity that $(M^t)_{ij}$ gives the probability of going from $S_i$ to $S_j$ in $t\geq 1$ steps, and this is denoted by
\(m_{ij}^{(t)}=\P[X_{t}=S_j~|~X_0=S_i].\)

Turning to probabilistic zero forcing, fix a connected graph $G$ and state space\footnote{In probabilistic zero forcing it is typical to omit the states which cannot be reached from some fixed initial blue vertex $\{v\}$, as well as collapsing states which behave analogously into a single state.} $\mathcal{S}=\{S_1,\ldots,S_s\}$.
Let the first state $S_1$ correspond to some fixed set of vertices $B$ being blue, the last state $S_s$ correspond to all vertices being blue, and the states $S_i$, $2\leq i\leq s-1$, correspond to some intermediate colorings.
Starting from state $S_1$, the probability that all vertices are blue at the end of round $r$ is $m_{1s}^{(r)}$. It follows that 
\[\P[\operatorname{pt}_{pzf}(G,B)=r]=m_{1s}^{(r)}-m_{1s}^{(r-1)}\]
and thus, as noted in \cite[Remark 2.12]{GH18}, the expected propagation time of $B$ is
\[\operatorname{ept}(G,B)=\sum_{r=1}^\infty r\left(m_{1s}^{(r)}-m_{1s}^{(r-1)}\right).\]
Letting $\mathbf{e}_s=[0,\ldots,0,1]^T$, this was later calculated explicitly in \cite[Theorem 2.2]{CC+20} to be
\[\operatorname{ept}(G,B)=((M-\mathds{1}\mathbf{e}_s^T-I)^{-1})_{1s}+1,\]
recalling $\mathds{1}=[1,\ldots,1]^T$.

\section{Reversion Probabilistic Zero Forcing}\label{sec:RPZF}
This section introduces reversion probabilistic zero forcing, a modification of the probabilistic zero forcing process where blue vertices have the chance to revert back to being white at the end of each round. 
Two variations of this process are defined. The first of which has a single stopping state when all vertices are white, and the second of which has two stopping states when the graph is either entirely white or entirely blue. 
The processes are then formulated as Markov chains before being contrasted with other similar growth/decay dynamics existing in the probability literature. 

\subsection{Single Absorption Reversion Probabilistic Zero Forcing}\label{sec:SARPZF}
Single absorption reversion probabilistic zero forcing (SARPZF) adds a second phase to each round of probabilistic zero forcing. We describe one round of SARPZF. 

\begin{definition}{\rm
Given a graph $G$ and set $B$ of currently blue vertices, in \defn{phase 1} each blue vertex $u\in B$ attempts to force each of its white neighbors $w\in N(u)\setminus B$ independently with probability \[\P (u\to w)=\frac{|N[u]\cap B|}{\deg{u}},\]
as in probabilistic zero forcing. We now have an updated set $B'$ of blue vertices. In \defn{phase 2}, each blue vertex $u\in B'$ \defn{reverts} (changes to being white) independently with probability $p\in(0,1)$. Phases 1 and 2, taken consecutively, define the \defn{single absorption reversion probabilistic zero forcing color change rule} (SARPZF color change rule). A round of \defn{single absorption reversion probabilistic zero forcing} (SARPZF) is one application of the SARPZF color change rule.}
\end{definition}

Unlike probabilistic zero forcing, $G$ is \textbf{not} guaranteed to be forced blue under SARPZF.
In fact, we will see in \cref{thm:absorbing} that, given enough time, SARPZF will always result in all vertices being white. We say SARPZF \defn{dies out} when this occurs.

\subsection{Dual Absorption Reversion Probabilistic Zero Forcing}\label{sec:DARPZF}
Notice that SARPZF may lead to $G$ being entirely blue any number of times before dying out. A natural question to ask is will $V(G)$ ever be entirely blue? If so, when is the first time we expect this to happen? To answer these questions, we define dual absorption reversion probabilistic zero forcing by introducing a stopping condition to SARPZF. 

\begin{definition}{\rm
The \defn{dual absorption reversion probabilistic zero forcing color change rule} (DARPZF color change rule) is defined by modifying the SARPZF color change rule as follows: after phase 1, if the set of currently blue vertices $B'$ is the entire vertex set, then no vertices revert in phase 2. A round of \defn{dual absorption reversion probabilistic zero forcing} (DARPZF) is one application of the DARPZF color change rule.}
\end{definition}
Collectively, SARPZF and DARPZF are referred to as \defn{reversion probabilistic zero forcing}, or RPZF.
We say that DARPZF \defn{fully forces} $G$ when every vertex is blue and \defn{dies out} when every vertex is white. We say DARPZF is \defn{absorbed} whenever it dies out or fully forces $G$.
In SARPZF, \defn{absorbed} refers only to SARPZF dying out, hence the terminology \textit{single} and \textit{dual} absorption.

\begin{remark}{\rm
We adopt the convention that $0<p<1$ in RPZF. 
When $p=1$, SARPZF trivially dies out in one step. 
When $p=0$, many of the results presented can be adapted to recover PZF results. However, most of these results are already known and so we refer the reader to the PZF literature.
Moreover, some results, such as those involving the matrix $Q_S$ introduced in \cref{eqn:partition}, are not suitable for adaptation to PZF.}
\end{remark}

\subsection{Markov chains for SARPZF and DARPZF}
As previously defined in the literature \cite{CC+20}, the \defn{simple state space of $G$} is the set of all white-blue colorings of $V(G)$, and we call each of these colorings a \defn{simple state}.
We combine simple states that behave analogously under RPZF into one \defn{state of $G$} and omit states that cannot occur. For example, on the complete graph $K_n$ we use $n+1$ states $S_0,\ldots,S_n$, where state $S_i$ is the condition of there being exactly $i$ blue vertices.
It is often helpful to think of a state as a set of blue vertices, with the remaining vertices being white.
We will use these two notions as convenient.
An \defn{ordered state space for $G$}, denoted $\mathcal{S}=(S_0,\ldots,S_s)$, is an ordered list of all states of $G$ where $S_0$ is the die-out state (all vertices are white), $S_s$ is the fully forced state (all vertices are blue), and the remaining states are chosen in some order.
Let $|S_i|$ denote the number of blue vertices in the coloring of state $S_i$. We say $\mathcal{S}$ is \defn{properly ordered} if $|S_i|\leq |S_{i+1}|$ for all $i\in\{0,\ldots,s-1\}$.

A graph $G$, a properly ordered state space $\mathcal{S}=\{S_0,\ldots,S_s\}$, and probability $p\in(0,1)$ then determine the SARPZF and DARPZF Markov transition matrices $M_S(G)$ and $M_D(G)$, respectively. We suppress the dependence on $G$ when $G$ is clear from context. For all $0\leq i,j\leq s$, $(M_S)_{ij}$ and $(M_D)_{ij}$ give the probability of going from state $S_i$ to state $S_j$ in one application of the SARPZF and DARPZF color change rule, respectively. We index these matrices starting from 0 to align with $\mathcal{S}$ and can partition them by
\begin{equation}\label{eqn:partition}
M_S=\left[\begin{array}{c|c}
1 &0\cdots0 \\
\hline
\mathbf{r} & Q_S 
\end{array}\right]
\quad\text{and}\quad
M_D=\left[\begin{array}{c|c|c}
1 &0\cdots0 & 0\\
\hline
\bo{a}_1 & Q_D & \bo{a}_2\\
\hline
0 &0\cdots 0 & 1
\end{array}\right]
\end{equation}
where the the $i^\text{th}$ row and column correspond to the state $S_i$. 
A state $S_i$ is \textit{absorbing} if the Markov chain cannot leave that state once entered.
Observe that $\bo{r}$ and $\bo{a}_1$ correspond to the absorbing state $S_0$ where all vertices are white and hence no vertex can be forced blue. Symmetrically, $\bo{a}_2$ corresponds to the absorbing state $S_s$ in DARPZF where every vertex is blue and so by definition of DARPZF no reversions can occur. 
Additionally, for any $t\in\naturals$,
\begin{equation}
(M_S)^t=\left[\begin{array}{c|c}
1 &0\cdots0 \\
\hline
* & (Q_S)^t 
\end{array}\right]
\quad\text{and}\quad
(M_D)^t=\left[\begin{array}{c|c|c}
1 &0\cdots0 & 0\\
\hline
* & (Q_D)^t & *\\
\hline
0 &0\cdots 0 & 1
\end{array}\right].
\end{equation}

\begin{remark}{\rm
The SARPZF Markov chain has only one absorbing state, $S_0$, as seen by the first row of $M_S$. The DARPZF chain, on the other hand, has two absorbing states $S_0$ and $S_s$ as seen by the first and last rows of $M_D$. }
\end{remark}

Note that $(M_S)^t$ and $(M_D)^t$ are right stochastic for all $t\in\naturals$ since they are Markov transition matrices. Moreover, the RPZF color change rules define Markov chains whose dynamics occur in two phases and thus can also be described as the product of two right stochastic matrices
\begin{equation}\label{eqn:matrixprod}
M_S=F R_S\quad\text{and}\quad M_D=F R_D
\end{equation}
where $F$ describes phase 1 (forcing, using the probabilistic zero forcing color change rule) and $R\in\{R_S,R_D\}$ describes phase 2 (reversion of blue vertices). We utilize this formulation in \cref{sec:KnRPZF}.

\begin{remark}{\rm
    The RPZF transition matrices are indexed from $0$, and submatrices preserve the indexing of their parent matrix. So, for example, $Q_S$ and $Q_D$ are indexed from 1.}
\end{remark}

To formulate the Markov chains themselves, fix a graph $G$, an initial set of blue vertices $B$, and a properly ordered state space $\mathcal{S}=(S_0,\ldots,S_s)$ on $G$. The \defn{SARPZF Markov chain} $(X_t^S)$ is a sequence of random variables $X_0^S,X_1^S,\ldots$ 
on $\mathcal{S}$ with one-step transition matrix $M_S$. Similarly, the \defn{DARPZF Markov chain} $(X_t^D)$ is a sequence of random variables $X_0^D,X_1^D,\ldots$
on $\mathcal{S}$ with one-step transition matrix $M_D$. That is, for any $t\in\naturals$ and any $i,j\in\{0,\ldots,s\}$,
\[\P[X^S_{t+1}=S_j~|~X_t^S=S_i]=(M_S)_{ij}\quad\text{and}\quad \P[X^D_{t+1}=S_j~|~X_t^D=S_i]=(M_D)_{ij}.\]

Just as probabilistic propagation time and expected propagation time were introduced for probabilistic zero forcing in \cite{GH18}, we introduce new parameters of interest for reversion probabilistic zero forcing. For the following definitions, $G$ is a graph with properly ordered state space $(S_0,\ldots,S_s)$, $B$ is a set of blue vertices in $G$, and $(X_t^S)$ and $(X_t^D)$ are SARPZF and DARPZF Markov chains on $G$ with reversion probability $p\in(0,1)$.
\begin{definition}{\rm
The \defn{probabilistic time of absorption for $B$ under SARPZF} (respectively DARPZF), denoted $\probabs_{S}(G;B,p)$ ($\probabs_D(G;B,p)$), is the first time at which every vertex turns white in SARPZF (all blue or all white in DARPZF). In probabilistic notation,
\[
\probabs_S(G;B,p)=\min\{t\geq 0:X_t^S=S_0\}
\quad\text{and}\quad
\probabs_D(G;B,p)=\min\{t\geq0:X_t^D\in\{S_0,S_s\}\}
\]
where we define $\min\varnothing=\infty$. }
\end{definition}
\begin{definition}{\rm
The \defn{expected time of absorption for $B$ under SARPZF} (DARPZF) is the expected value of the probabilistic time of absorption for $B$ under SARPZF (DARPZF), and is denoted by
\[
\expabs_S(G;B,p)=\E[\probabs_S(G;B,p)]
\quad\text{and}\quad
\expabs_D(G;B,p)=\E[\probabs_D(G;B,p)].
\]}
\end{definition}
\begin{definition}{\rm
The \defn{expected time of absorption under SARPZF} (DARPZF) of a connected graph $G$ is the minimum of the expected time of absorption for $B$ under SARPZF (DARPZF) over all one-vertex sets $B$ of $V(G)$ and is denoted by
\[
\expabs_S(G;p)=\min\{\expabs_{S}(G;\{v\},p):v\in V(G)\}
\quad\text{and}\quad
\expabs_D(G;p)=\min\{\expabs_{D}(G;\{v\},p):v\in V(G)\}.
\]}
\end{definition}
When considering a general RPZF Markov chain, we may omit the $S$ and $D$ subscripts and superscripts.
In light of the dual nature of DARPZF it is natural to ask, given a starting state (coloring) $S_i$, for what probability $p$ does $G$ have an equal chance of dying out or being fully forced?
\begin{definition}{\rm
    The \defn{critical reversion probability}, denoted $p_D(G,S_i)$, is the reversion probability such that the DARPZF Markov chain has equal probability of dying out or fully forcing when starting from state $S_i$.}
\end{definition}
We show in \cref{sec:gen_rpzf} that $p_D(G,S_i)$ exists for all connected graphs $G$ and all non-absorbing states $S_i$.

\subsection{Contact Processes}\label{sec:contact}
Reversion probabilistic zero forcing can be viewed as a discrete-time analog of what is known in the probability literature as the contact process, a basic model for population growth and the spread of infection.
Analogous to discrete time Markov chains, we say that $(\eta_t)$, $t\geq 0$, is a continuous-time Markov process on the state space $\mathcal{S}$ if for any $0\leq u_0<u_1<\cdots<u_k$ and any states $S_j,S_i,S_{i_{k-1}},\ldots,S_{i_{0}}\in\mathcal{S}$,
\[\P[\eta_{u_{k}+t}=S_j~|~\eta_{u_k}=S_i,\eta_{u_{k-1}}=S_{i_{k-1}},\ldots,\eta_{u_0}=S_{i_0}]=\P[\eta_t=S_j~|~\eta_0=S_i].\]
That is, given the present state, the past does not influence the future behavior of the process. 
The \emph{susceptible-infected-susceptible contact process} (SIS contact process) on a (possibly infinite) graph $G=(V,E)$ with infection parameter $\lambda\geq0$ is a continuous-time Markov process $(\eta_t)$ with state space $\{0,1\}^{V}$, where a state $\eta\in\{0,1\}^V$ is a configuration of zeros and ones on the graph.
At any time $t\geq 0$, each vertex has status either 0 (``healthy") or 1 (``infected").
The state of the entire system at time $t$ is then described by $\eta_t:V\to\{0,1\}$ where $\eta_t(v)$ is the status of vertex $v$ at time $t$.
Finally, letting $v$ be a vertex and $\eta$ a configuration, we say that the contact process $(\eta_t)$ evolves according to the following local transition rates: at vertex $v$,
\begin{align*}
    0\to1\quad&\text{at rate}\quad \lambda \sum_{w\in N(v)}\eta(w)\\
    1\to 0\quad&\text{at rate}\quad 1.
\end{align*}
Specifically, these are the rates for exponential random variables whose value corresponds to the waiting time until vertex $v$ changes status. 
Thus, infected vertices recover after some (exponential distribution) time with mean 1 independent of its neighbors, and healthy vertices become infected at a rate linearly proportional to its number of infected neighbors. 

A common topic in study of contact processes is the infection parameter $\lambda$ and its relation to the process surviving or dying out. The contact process is said to die out if \[\P[\eta_t\not\equiv 0~\forall t\geq 0]=0,\] 
where $\eta_t\not\equiv 0$ means there exists a $v\in V$ such that $\eta_t(v)\neq0$. 
Otherwise the process is said to survive.
It is well known that on finite graphs, no matter the initial configuration or infection parameter, the process dies out \cite{MVY13}. 
We prove the analogous result for SARPZF in \cref{thm:absorbing}. 

The discretized SIS contact process can be formally described as follows. 
Let $G$ be a (possibly infinite) graph. 
At each time $t\in \naturals$, every infected vertex $v\in V(G)$ infects each of its neighbors independently with probability $\beta$. 
Simultaneously, every infected vertex $v$ at time $t$ recovers with probability $p$. 
The exact meaning of ``simultaneously" differs depending on the particular discretization being considered. 
For instance, some models allow a vertex which recovers to be reinfected during the same time-step, whereas  
others assert that a recovered vertex must remain recovered.
SARPZF, on the other hand, is defined such that a vertex has a chance of being infected and then immediately recovering before having the chance to infect its neighbors. 

Over time, various formulations have been proposed to model the discrete-time SIS contact process. 
We describe a few of them and then contrast them to SARPZF. 
We consider only SARPZF because contact processes do not have the additional stopping condition that DARPZF has.
For any vertex $v\in V(G)$ and time $t\in\naturals$, define $\mf{p}_v(t)$ to be the probability that $v$ is infected (blue) at time $t$, and define $$\mf{q}_v(t+1)=\prod_{x\in N(v)}(1-\beta\mf{p}_x(t))$$ to be the probability that $v$ does not receive infection (is not forced) at time $t$. 
Finally, let $p$ be the probability that a vertex recovers (reverts to white). 
Then Wang et al. proposed the model \cite{wang2003epidemic}
\begin{equation}\label{wang_model}
1-\mf{p}_v(t+1)=(1-\mf{p}_v(t))\mf{q}_v(t+1)+p\mf{p}_v(t)\mf{q}_v(t+1)+\frac{1}{2}p\mf{p}_v(t)(1-\mf{q}_v(t+1).
\end{equation}
The first term is the probability of vertex $v$ entering time $t+1$ healthy and then not being infected, the second term is the probability of $v$ entering time $t+1$ infected, recovering, then not being reinfected, and the final term assumes that half of the time a vertex will undergo a ``curing event" after being reinfected.
Notice that this interpretation of ``simultaneous" has recovery occurring before infection.
This model is explored in more detail in \cite{chakrabarti2008epidemic}. 
Later models do away with the $1/2$ probability ``curing event" assumption, as well as stating the dynamics in terms of the probability of being infected, $\mf{p}_v$, instead of the probability of being healthy, $1-\mf{p}_v$. For instance, G{\'o}mez et al. introduced the model \cite{gomez2010discrete}
\begin{equation}\label{gomez_model}
\mf{p}_v(t+1)=(1-p)\mf{p}_v(t)+(1-\mf{q}_v(t+1))(1-\mf{p}_v(t))+p(1-\mf{q}_v(t+1))\mf{p}_v(t),
\end{equation}
accounting for the cases of an infected vertex failing to recover, a susceptible vertex being infected, and an infected vertex recovering then becoming reinfected. This formulation also makes the assumption that a vertex which recovers at time $t$ can immediately be reinfected at time $t$.
Contrast \cref{gomez_model} with the model presented by Ahn and Hassibi in \cite{ahn2013global} where
\begin{equation}\label{ahn_model}
\mf{p}_v(t+1)=(1-p)\mf{p}_v(t)+(1-\mf{p}_v(t))(1-\mf{q}_v(t+1)),
\end{equation}
doing away with the $p(1-\mf{q}_v(t+1))\mf{p}_v(t)$ term and so asserting that a vertex that recovers cannot be reinfected in the same time step.
This formulation seems most true to the notion of ``simultaneous" vertex infection and recovery since vertices can only undergo one status change each time-step.
\Cref{ahn_model} can be further simplified by truncating the terms of $\mf{q}_v(t+1)=\prod_{x\in N(v)} (1-\beta\mf{p}_x(t))$ with powers of $\beta$ greater than $1$, giving
\begin{equation}\label{pare_model}
    \mf{p}_v(t+1)=(1-p)\mf{p}_v(t)+(1-\mf{p}_v(t))\beta\sum_{x\in N(v)}\mf{p}_x(t).
\end{equation}
Note that this approximation is better for smaller values of $\beta$.
Par{\'e} et al. demonstrate in \cite{pare2018analysis} how \cref{pare_model} directly matches the model derived from applying Euler's method to the continuous-time mean field approximation for the SIS contact process, as well as providing analysis on the accuracy of \cref{ahn_model,pare_model}.
For a graph $G$ on $n$ vertices, these models can all be used to solve for $\mf{p}_v$ numerically, from which tests of accuracy are typically derived. 
A commonly considered parameter is the graph's expected infection density $\rho_t$ at time $t$. 
Given an infection rate $\beta$ and recovery rate $p$, this is computed as 
\[\rho_t=\frac{1}{n}\sum_{v\in V(G)}\mf{p}_v(t).\]

So far, the models described have all been homogeneous, meaning that the infection rate $\beta$ and recovery rate $p$ are constant for all vertices. SARPZF, on the other hand, is more akin to a heterogenous model, wherein the infection rate and recovery rate depend on the vertex. In fact, SARPZF is somewhere in the middle, with a constant recovery rate but variable infection rate. Looking to model SARPZF in the same way as the above models, we find
\begin{equation}\label{SAPZF_model}
    \mf{p}_v(t+1)=(1-p)\mf{p}_v(t)+(1-p)(1-\mf{p}_v(t))(1-\mf{q}_v(t+1)).
\end{equation}
The first term is the case of $v$ being infected after time $t$ and not reverting during time $t+1$, and the second term is vertex $v$ be healthy after time $t$, infected at time $t+1$, and not reverting at time $t+1$.
Notice that \cref{SAPZF_model} is most similar to \cref{ahn_model} but differs in two ways. First, the function $\mf{q}_v$ is different. Second, the $(1-\mf{p}_v(t))(1-\mf{q}_v(t+1))$ is multiplied by $1-p$ due to the fact that in SARPZF, a vertex has the chance to revert in the same time step it is infected. This interpretation of ``simultaneous" has recovery occurring after infection. 
One of the key differences between the models in the literature and SARPZF is the infection rate $\beta$ and, subsequently, the probability that vertex $v$ is not infected by a neighbor $\mf{q}_v$. Let $G$ be a graph on $n$ vertices and let $(X_t)$ be the SARPZF Markov chain on the properly ordered state space $\mathcal{S}=(S_0,\ldots,S_s)$ with reversion probability $p$. Notice that in SARPZF, if $B$ is our set of infected vertices, then vertex $v$ is ``infected" (forced) with probability
\[\P[B\to v]=1-\P[B\not\to v]=1-\prod_{x\in B\cap N(v)}\P[x\not\to v]=1-\prod_{x\in B\cap N(v)}\left(1-\frac{|B\cap N[x]|}{\deg x}\right).\]
Thus if $X_t$ is the set of blue vertices at time $t$ and $X_0=S_i$, 
\[\mf{q}_v(t+1)=\P[X_t\not\to v].\]
Applying the law of total probability over the state space $\mathcal{S}$, 
\begin{align*}
    \mf{q}_v(t+1)&=\sum_{j=0}^s\P[X_t\not\to v~|~X_t=S_j]\P[X_t=S_j]\\
    &=\sum_{j=0}^s\prod_{x\in S_j\cap N(v)}\left(1-\frac{|S_j\cap N[x]|}{\deg x}\right)
    \prod_{x\in S_j}\mf{p}_x(t)\prod_{x\in V(G)\setminus S_j}(1-\mf{p}_x(t)).
\end{align*}

Notice in this case that the rate of infection is exponentially proportional to the number of infected neighbors, and moreover is also dependent on the number of infected vertices at distance 2 from $v$.
Intuitively, this means that an infected vertex is \textit{more infectious} when it has more infected neighbors.

\section{RPZF on general graphs}\label{sec:gen_rpzf}
In order to utilize Markov chain theory in the analysis of RPZF, some more standard probability notation and results are established in \cref{prob_background}. This is followed by methods of calculating RPZF parameters using Markov matrices in \cref{MarkovTheoryResults} before a couple results on the expected number of blue vertices.

\subsection{Probabilistic Background}\label{prob_background}
Let $(X_t)$ be an RPZF Markov chain on a graph $G$ with reversion probability $p\in(0,1)$, properly ordered state space $(S_0,\ldots,S_s),$ and Markov transition matrix $M$. As shorthand notation, 
\(\P_i[A]=\P[A~|~X_0=S_i]\)
denotes the probability of event $A$ given the chain starts from state $S_i$, and \(\E_i[Y]=\E_i[Y~|~X_0=S_i]\)
denotes the expected value of random variable $Y$ given the chain starts from state $S_i$.
In this notation, the \defn{transition probability} $m_{ij}=\P_i[X_1=S_j]$ gives the probability of going from $S_i$ to $S_j$ in one step and is by definition the $(i,j)$th entry of the transition matrix $M$. 

Let $T(j)=\min\{t\geq 1:X_t=S_j\}$ be the \defn{time of first arrival} to $S_j$ (sometimes called the \defn{first return time}).
Notice that the starting state $X_0$ is not considered.
Define 
\[\rho_{ij}=\P_i[T(j)<\infty]=\P[X_t=S_j~\text{for some $t\geq 1$}~|~X_0=S_i]\]
to be the probability that the RPZF chain enters state $S_j$ after starting from state $S_i$. In this language, we can formally define the critical reversion probability $p_D(G,S_i)$ 
as the DARPZF reversion probability such that $\rho_{i0}=\rho_{is}=1/2$. In words, $p_D(G,S_i)$ is the reversion probability such that, starting from state $S_i$, we have equal probability to enter $S_0$ (all white) or $S_s$ (all blue) in DARPZF. We can also classify states of the Markov chain:
\begin{itemize}
    \item If $\rho_{ii}=1$ then state $S_i$ is said to be \defn{recurrent}.
    \item If $\rho_{ii}<1$ then state $S_i$ is said to be \defn{transient}.
    \item If $\P_i[X_1=S_i]=1$ then state $S_i$ is said to be \defn{absorbing}.
\end{itemize}
As we will show in \cref{lem:statetype}, state $S_0$ is absorbing and states $S_i$, $i\in\{1,\ldots,s-1\}$, are transient for SARPZF and DARPZF.

\begin{lemma}\label{lem:statetype}
Let $G$ be a connected graph and $(S_0,\ldots,S_s)$ be a properly ordered state space of $G$ for an RPZF chain with reversion probability $p\in(0,1)$. Then the state $S_0$ is absorbing, the state $S_i$ is transient for all $i\in\{1,\ldots,s-1\}$, and the state $S_s$ is absorbing in DARPZF and transient in SARPZF.
\end{lemma}
\begin{proof}
    Let $(X_t)$ be an RPZF chain on $G$. It is immediate from the Markov matrices (see \cref{eqn:partition}) that $S_0$ is an absorbing state for SARPZF and DARPZF. 
    Now let $i\in\{1,\ldots,s-1\}$ and consider 
    \[\rho_{ii}=\P_i[T(i)<\infty]=1-\P_i[T(i)=\infty].\] 
    To show $S_i$ is transient it suffices to show $\P_i[T(i)=\infty]>0$. 
    That is, it suffices to show there is a nonzero probability that the chain never returns to state $S_i$ after starting in state $S_i$.
    Suppose $X_0=S_i$ and that after phase 1 of the RPZF color change rule there are $b$ vertices colored blue. 
    Since $p>0$, the chain moves to state $S_0$ in phase 2 with probability $p^b>0$, whereafter it can never reach state $S_i$ again. The only exception is if, in DARPZF, all vertices are forced blue in phase 1. In this case, the chain remains in state $S_s$ and so never returns to $S_i$. 

    Finally, consider the state $S_s$ where all vertices are blue. In the DARPZF chain it is immediate from \cref{eqn:partition} that $S_s$ is absorbing. On the other hand, in the SARPZF chain, if $G$ has $n$ vertices then $$\P_s[T(s)=\infty]\geq\P_s[X_1=S_0]=p^n>0$$ 
    and so $S_s$ is transient in SARPZF.
\end{proof}

Notice that, as a consequence of \cref{lem:statetype}, every state is either transient or absorbing and the only recurrent states are the absorbing states. Additionally, the submatrices $Q_S$ and $Q_D$ from \cref{eqn:partition} correspond exactly to transient states. 
\begin{remark}{\rm
If $p=0$ then $S_s$ is actually an absorbing state for SARPZF, but we forbid this so that $Q_S$ always corresponds to transient states. }
\end{remark}
\subsection{RPZF Parameters and Expected Behavior}\label{MarkovTheoryResults}
For a real $n\times m$ matrix $A=[a_{ij}]$, define the matrix norm 
\[\norm{A}=\norm{A}_\infty=\max_{1\leq i\leq n}\sum_{j=1}^m |a_{ij}|\]
to be the maximum absolute row sum of the matrix. It is well known that $\norm{\cdot}$ is a matrix norm.
In particular, for $A,B\in\reals^{n\times m}$ we have $\norm{AB}\leq\norm{A}\norm{B}$. If $A$ is square and $\rho(A)$ denotes the spectral radius of $A$, then $\rho(A)\leq\norm{A}$. We refer the reader to \cite[Chapter 5.6]{HJ2012} for details on the matrix norm facts used in this section.

\begin{theorem}\label{thm:absorbing}
Let $G$ be a connected graph, and let $(X_t)$ be an RPZF Markov chain with transition matrix $M\in\{M_S,M_D\}$. 
Then with probability $1$, $(X_t)$ enters an absorbing state in finite time. 
That is, for any $B\subseteq V(G)$, $\P[\probabs(G,B)<\infty]=1$. Furthermore, $Q^t\to O$ as $t\to\infty$ for $Q\in\{Q_S,Q_D\}$.
\end{theorem}
\begin{proof}
Let $G$ be a connected graph, let $B\subseteq V(G)$, and let $(X_t)$ an RPZF chain on $G$ with properly ordered state space $\{S_0,\ldots,S_s\}$, transition matrix $M$, and reversion probability $p\in(0,1)$. It is known (see e.g. \cite[Ch 1.3]{DurrettEssentials}) that if $S_i$ is a transient state, then after some finite amount of time $(X_t)$ will never (re)visit $S_i$. It follows that after some finite amount of time $(X_t)$ moves into a non-transient state. By \cref{lem:statetype} we know any non-transient state is an absorbing state, and hence $\P[\probabs(G,B)<\infty]=1$.

To show $Q^t\to O$ as $t\to\infty$, we first claim $\norm{Q}<1$. 
To see this, let $\mathbf{r},\mathbf{a}_1,$ and $\mathbf{a}_2$ be as in \cref{eqn:partition}. 
Then $\mathbf{r}$ and $\mathbf{a}_1+\mathbf{a}_2$ are entrywise strictly positive. 
Indeed, for all $i=1,\ldots,s$ the entry $\mathbf{r}_i$ denotes the probability of the SARPZF chain moving from transient state $S_i$ to absorbing state $S_0$ in one time-step. 
Since $p>0$, this is always achievable by the event of every blue vertex reverting in Phase 2, giving $\mathbf{r}_i>0$. 
Now let $i=1,\ldots,s-1$. Then $(\mathbf{a}_1)_i$ denotes the probability of the DARPZF chain moving from transient state $S_i$ to absorbing state $S_0$ in one time-step. 
Since $p>0$, this is almost always achievable by every blue vertex reverting, giving $(\mathbf{a}_1)_i>0$. 
The exception is when the DARPZF chain forces every vertex blue after Phase 1. 
But in this case $(\mathbf{a}_2)_i>0$ and so $(\mathbf{a}_1+\mathbf{a}_2)_i>0$. 
Recall that $M$ is right stochastic and notice that the row sum for every row corresponding to a transient state contains a contribution from $\mathbf{r}$ or $\mathbf{a}_1+\mathbf{a}_2$. 
This implies every row sum of $Q$ must be strictly less than $1$ and thus $\norm{Q}<1$.
Using sub-multiplicativity,  $\norm{Q^t}\leq\norm{Q}^t\to 0$ as $t\to\infty$ which implies $Q^t\to O$ as $t\to\infty$. 
\end{proof}

The next results in this section are adapted from \cite{prob:2003}. In particular, \cref{lem:FundMatrix} comes from \cite[Theorem 11.4]{prob:2003}, \cref{thm:ExpectedAbsorption} follows \cite[Theorem 11.5]{prob:2003}, and \cref{thm:DieVsAbs} is a reworking of \cite[Theorem 11.6]{prob:2003}. 
We are interested in how long an RPZF chain is expected to stay in transient states. For a connected graph $G$ with properly ordered state space $(S_0,\ldots,S_s)$, let $N_S$ and $N_D$ be matrices such that $(N_S)_{ij}$ and $(N_D)_{ij}$ are the expected number of visits to transient state $S_j$ when starting from transient state $S_i$ in the SARPZF and DARPZF chains, respectively. In other words, 
\begin{equation}\label{eqn:N}
(N_S)_{ij}=\E_i[|\{t\geq 0:X^S_t=S_j\}|]\quad\text{and}\quad(N_D)_{ij}=\E_i[|\{t\geq 0:X^D_t=S_j\}|].
\end{equation}
These matrices exist and can be calculated from the Markov transition matrix. Recall that for an event $A$, $\mathbf{1}[A]=1$ if $A$ occurs and $\mathbf{1}[A]=0$ if $A$ does not occur. 

\begin{lemma}\label{lem:FundMatrix}
Let $(X_t)$ be an RPZF Markov chain on a connected graph $G$ with properly ordered state space $(S_0,\ldots,S_s)$, transition matrix $M\in\{M_S,M_D\}$ with corresponding transient state matrix $Q$, and reversion probability $p\in(0,1)$. Let $N$ be defined as in \cref{eqn:N} and correspond to the choice of $M$. Then $N$ exists and $N=\sum_{k=0}^\infty Q^k=(I-Q)^{-1}$.
\end{lemma}
\begin{proof}
It is well known that $\sum_{k=0}^\infty Q^k=(I-Q)^{-1}$ provided $\sum_{k=0}^\infty Q^k$ converges. This occurs if $\rho(Q)<1$, where $\rho(Q)$ denotes the spectral radius of $Q$. 
We showed $\norm{Q}<1$ in \cref{thm:absorbing} since $p>0$ and so $\rho(Q)\leq\norm{Q}<1$. 

Now let $N=[n_{ij}]$ be defined by \cref{eqn:N}. 
To show $N=\sum_{k=0}^\infty Q^k$, let $S_i$ and $S_j$ be two arbitrary transient states.
Then
\(\P_i\big[\mathbf{1}[X_k=S_j]=1\big]=(Q^k)_{ij}\)
and since $\mathbf{1}[X_k=S_j]$ is an indicator random variable, $\E_i\big[\mathbf{1}[X_k=S_j]\big]=(Q^k)_{ij}$. Therefore, the expected number of times the chain has been in state $S_j$ after $t$ rounds, starting from state $S_i$, is
\[\E_i\left[\sum_{k=0}^t \mathbf{1}[X_k=S_j]\right]
=\sum_{k=0}^t \E_i\big[\mathbf{1}[X_k=S_j]\big]
=\sum_{k=0}^t (Q^k)_{ij}.\]
It follows that the total number of times the chain is expected to be in state $S_j$, when starting from $S_i$, is
\[n_{ij}=\E_i\left[\sum_{k=0}^\infty \mathbf{1}[X_k=S_j]\right]=\sum_{k=0}^\infty \E_i\big[\mathbf{1}[X_k=S_j]\big]=\sum_{k=0}^\infty (Q^k)_{ij}\]
where the expected value can be passed under the limit by Fubini's theorem (see \cref{sec:appendix}) since 
\[\sum_{k=0}^\infty\E_i\big[|\mathbf{1}[X_k=S_j]|\big]=\sum_{k=0}^\infty (Q^k)_{ij}=(I-Q)^{-1}_{ij}<\infty.\]
Thus $N=\sum_{k=0}^\infty Q^k$. 
\end{proof}

Observe that summing across the $i$th row of $N$ gives the expected number of rounds the chain spends in transient states, having started from state $S_i$. Given that the chain starts at time $t=0$, this leads to the following theorem.

\begin{theorem}\label{thm:ExpectedAbsorption}
Consider an RPZF Markov chain with matrix $N$ as in \cref{eqn:N}. Let $t_i$ denote the expected time of absorption, given that the chain starts from transient state $S_i$, and let $\bo{t}$ be the column vector whose $i$th entry is $t_i$. Then $\bo{t}=N\mathds{1}$ where $\mathds{1}$ is the vector containing all ones.
\end{theorem}
Given a graph $G$ and properly ordered state space $(S_0,\ldots,S_s)$, one can use \cref{thm:ExpectedAbsorption} to calculate $\expabs(G;S_i,p)=\bo{t}_i$. 
This tells us how long we expect to wait until the RPZF chain enters an absorbing state. 
But in the case of DARPZF, is that absorbing state more likely to be all white or all blue? 
This is answered by the next result.

\begin{theorem}\label{thm:DieVsAbs}
Consider the DARPZF Markov chain on a connected graph $G$ with properly ordered state space $(S_0,\ldots,S_s)$ and transition matrix $M_D$. 
Let $C=[c_{ij}]$ be the $(s-1)\times 2$ matrix such that, starting the chain from transient state $S_i$, $c_{i1}$ is the probability that the graph dies out and $c_{i2}$ is the probability that the graph is fully forced. Then
\[C=(I-Q_D)^{-1}[\bo{a}_1~\bo{a}_2]\]
where the matrix $Q_D$ and vectors $\bo{a}_1$, $\bo{a}_2$ come from \cref{eqn:partition}.
\end{theorem}
\begin{proof}
Let $A=[\bo{a}_1~\bo{a}_2]$ where $A$ inherits the row indices of $M_D$ and let $m_{ij}=(M_D)_{ij}$ denote the probability of going from state $S_i$ to state $S_j$ in one step of the DARPZF Markov chain. 
Then by definition of $C$,
\begin{equation}\label{eqn:AbsProb1}
c_{i1}=m_{i0}+\sum_{k=1}^{s-1}m_{ik}c_{k1}.
\end{equation}
Indeed, starting from state $S_i$, the chain can either die out in one step with probability $m_{i0}$ or move into to some transient $S_k$ and eventually die out with probability $m_{ik}c_{k1}$. Similarly, 
\begin{equation}\label{eqn:AbsProb2}
c_{i2}=m_{is}+\sum_{k=1}^{s-1}m_{ik}c_{k2}.
\end{equation}
Observe that $m_{ik}=(Q_D)_{ik}$ for $1\leq i,k\leq s-1$.
Moreover, writing $A=[a_{ij}]$, notice from \cref{eqn:partition} that $m_{i0}=a_{i1}$ and $m_{is}=a_{i2}$. Substituting all this into \cref{eqn:AbsProb1,eqn:AbsProb2},
\begin{equation*}
    c_{i1}=a_{i1}+\sum_{k=1}^{s-1}(Q_D)_{ik}c_{k1}=a_{i1}+(Q_D C)_{i1}
\end{equation*}
and
\begin{equation*}
    c_{i2}=a_{i2}+\sum_{k=1}^{s-1}(Q_D)_{ik}c_{k2}=a_{i2}+(Q_D C)_{i2}.
\end{equation*}
It follows that $C=A+Q_D C$ from which we conclude $C=(I-Q_D)^{-1}A$ as claimed, where $(I-Q_D)^{-1}$ exists by \cref{lem:FundMatrix}.
\end{proof}
\begin{corollary}
For any connected graph $G$ with properly ordered state space $(S_0,\ldots, S_s)$ and for any transient state $S_i$, the critical reversion probability $p_D(G,S_i)$ exists. That is, there exists a reversion probability such that $\P_i[T(0)<\infty]=\P_i[T(s)<\infty]=1/2$, where $T(j)=\inf\{t\geq 1:X_t^D=S_j\}.$
\end{corollary}
\begin{proof}
    Let $C=[c_{ij}]$ be as in \cref{thm:DieVsAbs}. The entries of $Q_D$ are continuous functions in $p$ and $I-Q$ is invertible on $p\in(0,1)$. Hence the entries of $C=(I-Q)^{-1}$ are continuous functions in $p\in(0,1)$. By taking $p$ sufficiently small, $c_{i1}=c_{i1}(p)<1/2$. Similarly, by taking $p$ sufficiently large, $c_{i1}=c_{i1}(p)>1/2.$ By the intermediate value theorem, $p_D(G,S_i)$ exists so that $c_{i1}=c_{i1}(p_D(G,S_i))=1/2$. Since $c_{i1}=\P_i[T(0)<\infty]$, $c_{i2}=\P_i[T(s)<\infty]$, and $c_{i2}=1-c_{i1}$, the result follows.
\end{proof}

We now follow the methods of Theorem 3.1 in \cite{CC+20} to calculate the expected number of blue vertices after one step of the SARPZF Markov chain. Recall that a state $S_i$, and thus a random variable $X_t$ of an RPZF chain, can be viewed as a set of blue vertices.

\begin{proposition}\label{OneStepExpected}
    Let $(X_t)$ be a SARPZF Markov chain on a connected graph $G$ with reversion probability $p\in(0,1)$ and properly ordered state space $(S_0,S_1,\ldots,S_s)$. Let $F_t$ be the number of vertices forced during phase 1 of time $t$ (before reversion). Then for all $i\in\{0,\ldots,s\}$, $\E_i[|X_1|]=(1-p)(|S_i|+\E_i[F_1])$. 
\end{proposition}
\begin{proof}
    Suppose $X_0=S_i$. Let $v_1,\ldots,v_{|S_i|+F_1}$ be the vertices that are blue after phase 1 of the SARPZF color change rule at time $t=1$. 
    Notice that $\mathbf{1}[v_j\not\in X_1]$ are i.i.d. indicator random variables for the event of $v_j$ reverting, $j=1,\ldots,|S_i|+F_1$. Thus
    \begin{align*}
        \E_i[|X_1|]&=|S_i|+\E_i F_1-\E_i\left[\sum_{j=1}^{|S_i|+F_1} \mathbf{1}[v_j\not\in X_1]\right]\\
        &=|S_i|+\E_i F_1-\E_i\big[(|S_i|+F_1)\mathbf{1}[v_1\not\in X_1]\big]\\
        &=|S_i|+\E_i F_1-\E_i[|S_i|+F_1]p\\
        &=(1-p)(|S_i|+\E_i F_1)
    \end{align*}
    since $v_1$ reverts with probability $p$.
\end{proof}

\begin{proposition}\label{lem:OneStepUpperbound}
    Let $G$ be a connected graph on $n$ vertices. Let $B=\{v_1,v_2,\ldots,v_b\}$ be the set of blue vertices with $b\geq 1$ and let $p\in(0,1)$ be the reversion probability. The expected number of blue vertices in $G$ after one step of the SARPZF color change rule is bounded above by $(1-p)(b+b^2)$.
\end{proposition}
\begin{proof}
    Let $(X_t)$ be the SARPZF chain on $G$ with starting state $X_0=S_i$ corresponding to the vertices $B=\{v_1,\ldots,v_b\}$ colored blue. Let $F_t$ be the number of vertices forced during phase 1 of time $t$ (before reversion).
    We start by bounding $\E_i F_1$.
    For each $j=1,\ldots,b$, let $B_j$ be the set of vertices forced by $v_j$ in phase 1. 
    Note that multiple blue vertices may force the same white vertex and so the $B_j$'s may intersect. Then 
    \[\E_i F_1= \E_i[|B_1\cup\cdots\cup B_b|]\leq \E_i\left[\sum_{j=1}^b |B_j|\right]=\sum_{j=1}^b\E_i[|B_j|].\]
    Calculating this expected value for each $j$,
    \begin{align*}
    \E_i[|B_j|]&=\sum_{k=0}^{|N(v_j)\setminus B|}k\P[|B_j|=k]\\
    &=\sum_{k=0}^{|N(v_j)\setminus B|}k\binom{|N(v_j)\setminus B|}{k}\left(\frac{N[v_i]\cap B|}{\deg(v_j)}\right)^k\left(1-\frac{N[v_j]\cap B|}{\deg(v_j)}\right)^{|N(v_j)\setminus B|-k}
    \end{align*}
    which is the expected value of a binomially distributed random variable. Thus
    \[\E_i[|B_j|]=|N(v_j)\setminus B|\left(\frac{N[v_j]\cap B|}{\deg(v_j)}\right)\leq |N[v_j]\cap B|\]
    and we conclude
    \(\E_i F_1\leq \sum_{j=1}^b|N[v_j]\cap B|\leq b|B|=b^2.\) It now follows from \cref{OneStepExpected} that 
    $$\E_i[|X_1|]=(1-p)(b+b^2)$$ 
    since $|S_i|=b$.
\end{proof}

\begin{remark}{\rm
    For any connected graph with SARPZF and DARPZF Markov processes $(X_t^S)$ and $(X_t^D)$, $\E_i[|X_1^D|]\geq\E_i[|X_1^S|]$ for all states $S_i$. This can be seen by recalling from \cref{eqn:matrixprod} that SARPZF and DARPZF share the same distribution during phase 1. In phase 2 of SARPZF, each vertex reverts with probability $p$. In phase 2 of DARPZF each vertex reverts with probability $p$ unless the chain steps into the absorbing state $S_s$, in which case each vertex reverts with probability 0. Hence SARPZF is expected to have at least as many reversions as DARPZF. }
\end{remark}

\section{Asymptotic thresholds for RPZF}\label{sec:thresholds}
In this section we consider the behavior of RPZF as the number of vertices grows towards infinity. In particular, we deduce thresholds for the number of blue vertices required to fully force the complete, complete bipartite, and star graphs in one step with high probability.
When utilizing the RPZF Markov chain $(X_t)$ on these graphs, it will be convenient to consider the random variables $X_t$ as taking integer values representing the number of blue vertices in the graph at the end of time $t$. 

\subsection{The complete graph}\label{sec:KnRPZF}
In any RPZF chain, the outcomes of phase 1 (forcing) and phase 2 (reversion) depend only on what vertices are blue at the start of the phase. 
It follows that we can decompose $M\in\{M_S,M_D\}$ as the product of two stochastic matrices $M=FR$ where $R\in\{R_S,R_D\}$ corresponds to the choice of $M$. 
The explicit formulas of $F=[f_{ij}]$ and $R=[r_{ij}]$ for the complete graph on $n$ vertices, denoted $K_n$, are given in \cref{thm:KnDef}. 
Before that, we first provide the known Markov transition matrix for (traditional) probabilistic zero forcing on $K_n$.

\begin{theorem}{\rm \cite[Theorem 2.4]{CC+20}\label{thm:KnPZF}}
    Let $(S_1,\ldots,S_n)$ be the properly ordered state space where $S_k$ is the state of having $k$ blue vertices in $K_n$, $n\geq 2$. The $n\times n$ Markov transition matrix $K(n)=[k_{ij}]$ for probabilistic zero forcing on $K_n$ is given by
    \[
k_{ij}=
\begin{cases}
\binom{n-i}{j-i}\left(1-\left(1-\frac{i}{n-1}\right)^i\right)^{j-i}\left(\left(1-\frac{i}{n-1}\right)^i\right)^{n-j}, &1\leq i\leq j\leq n\\
0, &\text{otherwise}
\end{cases}
\]
where $0^0=0!=1$.
\end{theorem}

The idea behind \cref{thm:KnPZF} is as follows. Let $B$ be the set of currently blue vertices on $K_n$ with $|B|=b$. Then for any $v\in B$ and $w\not\in B$, $\P[v\rightarrow w]=\frac{b}{n-1}$ and $\P[v\not\rightarrow w]=1-\frac{b}{n-1}$. Now for any given $w\not\in B$, each $v\in B$ will independently attempt to force $w$. Hence
\begin{equation}
\P[\forall v\in B,v\not\rightarrow w]=\left(1-\frac{b}{n-1}\right)^b\quad\text{and}\quad\P[B\rightarrow w]=1-\left(1-\frac{b}{n-1}\right)^b,\label{eqn:ForcingProb}
\end{equation}
where $B\rightarrow w$ denotes the event that $w$ is forced by some $v\in B$. 
For notational convenience, we denote the probability that $b$ blue vertices force a vertex in $K_n$ by
\begin{equation}\label{eqn:qnb}
q(n,b)=1-\left(1-\frac{b}{n-1}\right)^b
\end{equation}

\begin{observation}{\rm
With this notation, the matrix $K(n)=[k_{ij}]$ is given by 
$$k_{ij}=\binom{n-i}{j-i}\left(q(n,i)\right)^{j-i}\left(1-q(n,i)\right)^{n-j}$$ when $1\leq i\leq j$ and $0$ otherwise.}
\end{observation}
 We now define the transition matrices for SARPZF and DARPZF using \cref{thm:KnPZF}. For two matrices $A$ and $B$, the direct sum $A\oplus B$ is defined by
 \[A\oplus B=\begin{bmatrix}A&O\\O&B\end{bmatrix}.\]

\begin{theorem}\label{thm:KnDef}
Let $(S_0,\ldots,S_n)$ be the properly ordered state space where $S_k$ is the state of having $k$ blue vertices in $K_n$, $n\geq 2$, and let $p\in(0,1)$. 
Then we can write the $(n+1)\times(n+1)$ SARPZF and DARPZF Markov transition matrices for $K_n$ as $M_S=\left([1]\oplus K(n)\right)R_S$ and $M_D=\left([1]\oplus K(n)\right)R_D$ where $\left([1]\oplus K(n)\right)$, $R_S$, and $R_D$ are $(n+1)\times (n+1)$ right stochastic matrices. 
The matrices $R_S$ and $R_D$ describe reversion in SARPZF and DARPZF respectively and are given by
\[(R_S)_{ij}=\begin{cases}
\binom{i}{j}p^{i-j}(1-p)^j,& 0\leq j\leq i\leq n\\
0,& \text{otherwise}
\end{cases}
\]
and
\[
(R_D)_{ij}=\begin{cases}
\binom{i}{j}p^{i-j}(1-p)^j, & 0\leq i<n\text{ and }0\leq j\leq i\\
1, &i=j=n\\
0, &\text{otherwise.}\\
\end{cases}
\]
\end{theorem}
\begin{proof}
    The matrix $K(n)$ comes from \cite[Theorem 2.4]{CC+20} and describes probabilistic zero forcing for the complete graph on $n$ vertices, where the first row and column correspond to the graph having one blue vertex. Hence $F=[1]\oplus K(n)$ so that the first row and column of $F$ correspond to the absorbing state where the graph has zero blue vertices.

    For $0\leq j\leq i\leq n$, we define $(R_S)_{ij}$ and $(R_D)_{ij}$ to be the probability of reverting from $i$ blue vertices to $j$ blue vertices in phase 2 of SARPZF and DARPZF, respectively; that is, the probability that $i-j$ blue vertices turn white. Given a set of $i$ blue vertices, the probability that a particular collection of $i-j$ vertices revert and the remaining $j$ do not is $p^{i-j}(1-p)^j$, and there are $\binom{i}{j}$ ways to pick the $j$ vertices which do not revert. When $i=j=n$ in DARPZF, the probability of no vertices reverting is $1$.
\end{proof}

Note that $K(n)$ is indexed from 1, whereas $M_S$, $M_D$, $F$, $R_S$, and $R_D$ are indexed from 0. We say that the $i$th row/column corresponds to the state $S_i$, which in turn corresponds to $K_n$ having exactly $i$ blue vertices.
We next provide an explicit description of the one-step Markov transition probability $\P_b[X_1=k]$. To state the result, recall the Kronecker delta where $\delta_{ij}=1$ if $i=j$ and $0$ if $i\neq j$.
\begin{proposition}\label{PMF}
    Let $(X_t^S)$ and $(X_t^D)$ be the SARPZF and DARPZF Markov chains on $K_n$, $n\geq 2$, with reversion probability $p\in(0,1)$,
    and let $q(n,b)$ be defined as in \cref{eqn:qnb}. Then for any $1\leq b\leq n$ and $0\leq k\leq n$, $\P_b[X_1^S=k]$ is equal to both of the following:
    \begin{enumerate}
        \item[(1)] \(\displaystyle\sum_{i=\max\{b,k\}}^n\binom{n-b}{i-b}\binom{i}{k}(1-p)^k p^{i-k}q(n,b)^{i-b}(1-q(n,b))^{n-i}\)
        \item[(2)] \(\displaystyle\sum_{i=0}^{\min\{b,k\}} \binom{n-b}{k-i}\binom{b}{i}(1-p)^ip^{b-i}[(1-p)q(n,b)]^{k-i}[1-(1-p)q(n,b)]^{n-b-(k-i)}.\)
    \end{enumerate}
    Additionally,
    \begin{align*}
    \P_b[X_1^D=k]
    &=\P_b[X_1^S=k]-\binom{n}{k}p^{n-k}(1-p)^k q(n,b)^{n-b}+\delta_{nk}\left(q(n,b)\right)^{n-b}
    .\end{align*}
\end{proposition}

{\em Proof.}
    Let $M_S=FR_S$ and $M_D=FR_D$ denote the Markov transition matrices for the SARPZF and DARPZF chains $(X_t^S)$ and $(X_t^D)$. Then
    \[\P_b[X_1^S=k]=(M_S)_{bk}=\sum_{i=0}^n F_{bi}(R_S)_{ik}=\sum_{i=\max\{b,k\}}^n F_{bi}(R_S)_{ik}\]
    since $F_{bi}=0$ when $i<b$, and $(R_S)_{ik}=0$ when $i<k$.
    Observing $b\geq 1$, substitute
    $$F_{bi}=\binom{n-b}{i-b}\left(q(n,b)\right)^{i-b}\left(1-q(n,b)\right)^{n-i}$$ and 
    $$(R_S)_{ik}=\binom{i}{k}p^{i-k}(1-p)^k$$
    to get formula (1).
    
    Considering $\P_b[X_1^D=k]$, again $F_{bi}=0$ when $i<b$ and $(R_D)_{ik}=0$ when $i<k$. Hence
    \[\P_b[X_1^D=k]=\sum_{i=\max\{b,k\}}^{n} F_{bi}(R_D)_{ik}\]
    and observing that $(R_S)_{ik}=(R_D)_{ik}$ for $k\leq i\leq n-1$ gives
    \[\P_b[X_1^D=k]=\P_b[X_1^S=k]-F_{bn}(R_S)_{nk}+F_{bn}(R_D)_{nk}.\]
    Using the fact that $(R_D)_{nk}=0$ when $k<n$ and $(R_D)_{nn}=1$, this simplifies to
    \[\P_b[X_1^D=k]=\P_b[X_1^S=k]-\binom{n}{k}p^{n-k}(1-p)^k q(n,b)^{n-b}+\delta_{nk}\left(q(n,b)\right)^{n-b}.\]

    To derive formula (2) for $\P_b[X_1^S=k]$, consider $B$, the set of currently blue vertices with $|B|=b$, and $V\setminus B$, the set of currently white vertices.
    Each vertex $w\in V\setminus B$ is blue after one application of the SARPZF color change rule with probability $(1-p)q(n,b)$ since it must be forced blue with probability $q(n,b)$ and not revert with probability $1-p$.
    Each vertex $v\in B$ is blue after one application of the SARPZF color change rule with probability $1-p$.
    Suppose $X_0^S=b$ and $X_1^S=k$. Let $i$ denote the number of vertices in $B$ that did not revert at time $t=1$.
    There are $\binom{b}{i}$ choices of vertices to not revert, each case occurring with probability $(1-p)^i p^{b-i}$.
    Then $k-i$ vertices in $V\setminus B$ were forced blue and did not revert.
    This can happen in one of $\binom{n-b}{k-i}$ ways, each with probability $[(1-p)q(n,b)]^{k-i}[1-(1-p)q(n,b)]^{n-b-(k-i)}.$ 
    Summing over all choices for $i$, $\P_b[X_1^S=k]$ is equal to
    \[\sum_{i=0}^{\min\{b,k\}} \binom{n-b}{k-i}\binom{b}{i}(1-p)^ip^{b-i}[(1-p)q(n,b)]^{k-i}[1-(1-p)q(n,b)]^{n-b-(k-i)}. ~~~~~~ \cvd\] 

Note that SARPZF can also be described in terms of known probability distributions. 
To see how, let $B$ be the set of currently blue vertices with $|B|=b$, let $X$ be a random variable equal to the number of vertices from $B$ which do not revert, and let $Y$ be a random variable equal to the number of vertices from $V\setminus B$ which are forced blue and do not revert.
Then $X$ and $Y$ are independent with $X\sim \operatorname{Binomial}(b,1-p)$ and $Y\sim \operatorname{Binomial}(n-b,(1-p)q(n,b)).$ Moreover, if $X_0^S=b$ then $X_1^S=X+Y$ follows a Poisson binomial distribution with $p_1=\cdots=p_b=1-p$ and $p_{b+1}=\cdots=p_n=(1-p)q(n,b)$. The probability mass function for this distribution simplifies to formula (2) in \cref{PMF}.

\begin{theorem}\label{OneStepSAR}
Suppose $K_n$ has $1\leq b\leq n-2$ vertices colored blue and RPZF chain $(X_t)$ with reversion probability $p\in(0,1)$. Then the probability $K_n$ dies out in one step of the RPZF chain converges to $p^b e^{b^2(p-1)}$ as $n\to\infty$.
\end{theorem}
\begin{proof}
By formula (2) of \cref{PMF}, 
\[\P_b[X_1^S=0]=p^b[1-(1-p)q(n,b)]^{n-b}=p^b\left[p+(1-p)\left(1-\frac{b}{n-1}\right)^b\right]^{n-b}.\]
Similarly, 
\[\P_b[X_1^D=0]=p^b\left[p+(1-p)\left(1-\frac{b}{n-1}\right)^b\right]^{n-b}-p^n\left(1-\left(1-\frac{b}{n-1}\right)^b\right)^{n-b}.\]
To calculate the limits of these values, notice first that
\(0\leq\left(1-\left(1-\frac{b}{n-1}\right)^b\right)^{n-b}\leq1\)
for all $0\leq b\leq n-2$ and hence
\[0\leq p^n\left(1-\left(1-\frac{b}{n-1}\right)^b\right)^{n-b}\leq p^n\to0\]
as $n\to\infty$ since $p\in(0,1)$. 
We are left to consider $p^b\left(p+(1-p)\left(1-\frac{b}{n-1}\right)^b\right)^{n-b}.$ 
Taking the limit as $n\to\infty$, this is equal to
\[p^b\operatorname{exp}\left\{\lim_{n\to\infty} \frac{1}{(n-b)^{-1}} \log\left[p+(1-p)\Big(1-\frac{b}{n-1}\Big)^b\right]\right\}\]
and applying L'Hôpital's rule this simplifies to $p^be^{b^2(p-1)}.$
\end{proof}

\begin{remark}{\rm
Let $(X_t^D)$ be the DARPZF Markov chain on $K_n$. The probability of $K_n$ being fully forced in DARPZF when starting from $b$ blue vertices is bounded above by $1-\P_b[X_1^D=0]$. Indeed, 
\[\P_b[\min\{t:X_t^D=n\}<\infty]=1-\P_b[\min\{t:X_t^D=0\}<\infty]\leq 1-\P_b[X_1^D=0].\]}
\end{remark}

The remainder of this section is dedicated to threshold-like results for RPZF. These results concern the necessary number of blue vertices $b_n$ for a particular event to occur in one step of the RPZF chain with high probability, where $b_n$ is a function of $n$, the total number of vertices in the graph, and $n\to\infty$. We first consider the expected number of blue vertices after one step of an RPZF process on $K_n$, starting from $b$ blue vertices.
\begin{theorem}\label{OneStepExact}
    Let $(X_t^S)$ and $(X_t^D)$ be the SARPZF and DARPZF Markov chains on $K_n$ with reversion probability $p\in(0,1)$, and let $q(n,b)$ be defined as in \cref{eqn:qnb}. Then
    \[\E_b[X_{1}^S]=(1-p)\left(b+(n-b)q(n,b)\right)\]
    and
    \[\E_b[X_1^D]=\E_b[X_{1}^S]+np\,q(n,b)^{n-b}.\]
\end{theorem}
\begin{proof}
    We consider first the SARPZF Markov chain. Suppose that at time $0$ we have $b$ blue vertices and $w=n-b$ white vertices. 
    For each white vertex $v_1,\ldots,v_w$, 
    let $\mathbf{1}[v_i]=1$ if $v_i$ is colored blue in phase 1 of the current round and 0 otherwise, and define $F_1=\sum_{i=1}^w \mathbf{1}[v_i]$.
    Then $F_1$ is the number of vertices forced blue at time $t=1$, and since the $\mathbf{1}[v_i]$'s are i.i.d we have
    \[\E_b[F_1]=w\E_b[\mathbf{1}[v_i]]=(n-b)\left(1-\left(1-\frac{b}{n-1}\right)^{b}\right)=(n-b)q(n,b)\label{eqn:F0}\]
    where $\E[\mathbf{1}[v_i]]$ comes from \cref{eqn:ForcingProb}. The result now follows by substituting into $\E_b[X_1^S]=(1-p)(b+\E_b[F_1])$ from \cref{OneStepExpected}.

    To calculate $\E_b[X_1^D]$, first observe that by \cref{OneStepExpected} and linearity,
    \[\E_b[X_1^S]=(1-p)\E_b[F_1+b]=(1-p)\sum_{k=b}^n k\P_b[F_1+b=k].\]
    By definition of DARPZF there is no reversion when $F_1+b=n$. Hence
    \begin{align*}
    \E_b[X_1^D]&=(1-p)\sum_{k=b}^{n-1} k\P_b[F_1+b=k]+n\P_b[F_1+b=n]\\
    &=\E_b[X_1^S]-(1-p)n\P_b[F_1+b=n]+n\P_b[F_1+b=n]\\
    &=\E_b[X_1^S]+np\P_b[F_1+b=n]
    \end{align*}
    and substituting $\P_b[F_1=n-b]=q(n,b)^{n-b}$ finishes the proof.
\end{proof}

The next result gives a threshold number of blue vertices for the complete graph to completely forced in one step of RPZF with high probability. 
In particular, if $K_n$ has $b_n=\Omega(\sqrt{n\log n})$ blue vertices with constant $C>\sqrt{2}$ then $K_n$ is expected to have $n$ blue vertices after one application of the DARPZF color change rule. We refer the reader to \cref{sec:appendix} for a review of asymptotic notation and their common properties.
\begin{lemma}\label{OneStepThreshold}
    Let $(X_t^D)$ be the DARPZF Markov chain on $K_n$, $n\geq 2$, with reversion probability $p\in(0,1)$. If $b_n\geq\sqrt{n\log n^{2+\gamma}}$ with $\gamma>0$, then 
    \[\lim_{n\to\infty}|n-\E_b[X_1^D]|=0.\]
\end{lemma}

\begin{proof}
    Observe $\E_{n-1}[X_1^D]=\E_{n}[X_1^D]=n$, so we may assume $b_n\leq n-2$. Let $b_n\geq \sqrt{n\log n^{2+\gamma}}$ with $\gamma>0$. 
    Throughout this proof, let 
    \[g(n,b_n)=1-q(n,b_n)=\left(1-\frac{b_n}{n-1}\right)^{b_n}.\]
    Using \cref{OneStepExact}, observe that 
    \[\E_{b_n}[X_1^D]=(1-p)\big(b_n+(n-b_n)\left(1-g(n,b_n)\right)\big)+np\big(1-g(n,b_n)\big)^{n-b_n}
    \]
    which, after distributing $(n-b_n)$ and canceling the $b_n$ terms, is equal to
    \[(1-p)\big(n-(n-b_n)g(n,b_n)\big)+np\big(1-g(n,b_n)\big)^{n-b_n}.\]
    Now distribute the $(1-p)$ and simplify to get
    \begin{align*}
    \E_{b_n}[X_1^D]&=(1-p)n-(1-p)(n-b_n)g(n,b_n)+np\big(1-g(n,b_n)\big)^{n-b_n}\\
    &=n-(1-p)(n-b_n)g(n,b_n)-np\Big(1-\big(1-g(n,b_n)\big)^{n-b_n}\Big)
    \end{align*}
    Additionally, $\E_{b_n}[X_1^D]\leq n$ because $X_1^D\leq n$ and so $|n-\E_{b_n}[X_1^D]|=n-\E_{b_n}[X_1^D]$, which in turn simplifies to
    \[
    (1-p)(n-b_n)g(n,b_n)+np\Big[1-\big(1-g(n,b_n)\big)^{n-b_n}\Big].
    \]
    We show that each of these two terms converge to 0 as $n\to\infty$.
    
    For the first term it suffices to show that 
    \(
    ng(n,b_n)\to 0
    \)
    as $n\to\infty$ since $n-b_n$ and $g(n,b_n)$ are non-negative. 
    Recalling the Taylor expansion $\log(1-x)=-\sum_{k=1}^\infty x^k/k$ for $|x|<1$ we have
    \begin{align}\label{eqn:gnb}
    g(n,b_n)=\exp{b_n\log\left(1-\frac{b_n}{n-1}\right)}
    &=\exp{b_n\left(-\frac{b_n}{n-1}-\sum_{k=2}^\infty\frac{b_n^k}{k(n-1)^k}\right)}\nonumber\\
    &=\exp{-\frac{b_n^2}{n-1}}\exp{-\sum_{k=2}^\infty\frac{b_n^{k+1}}{k(n-1)^k}}
    \end{align}
    for all $b_n<n-1.$
    Then $g(n,b_n)\leq e^{-b_n^2/(n-1)}$ because $b_n>0$.
    Using the inequality $b_n\geq\sqrt{n\log n^{2+\gamma}}$, 
    \begin{equation}\label{eqn:n-delta}
    g(n,b_n)
    \leq\exp{-\frac{b_n^2}{n-1}}
    \leq\exp{-\frac{n\log n^{2+\gamma}}{n-1}}
    <\exp{-\frac{n\log n^{2+\gamma}}{n}}
    =n^{-(2+\gamma)}.
    \end{equation}
    Hence
    \[
        0\leq(1-p)\left(n-b_n\right)g(n,b_n)\leq \exp{-\frac{b_n^2}{n-1}}n<n^{-(1+\gamma)}\to 0
    \]
    as $n\to\infty$.
    
    It is left to show 
    \[\lim_{n\to\infty} np\Big(1-\big(1-g(n,b_n)\big)^{n-b_n}\Big)=0.\]
    Define $H(n)=(n-b_n)\log\left(1-g(n,b_n)\right)$ so that 
    \(1-\big(1-g(n,b_n)\big)^{n-b_n}=1-e^{H(n)}.\)
    It thus suffices to show that $n(1-e^{H(n)})\to 0$ as $n\to\infty$.
    Observe
    \[H(n)=(n-b_n)\log(1-g(n,b_n))\leq0\]
    since $0\leq g(n,b_n)<1.$
    Hence, to prove $n(1-e^{H(n)})\to 0$ it is enough to show 
    \[H(n)\geq\log\left(1-n^{-(1+\gamma)}\right)\]
    for sufficiently large $n$, since then
    \[0\leq n(1-e^{H(n)})\leq n(1-(1-n^{-(1+\gamma)}))=n^{-\gamma}.\]
    To this end, notice that because $H(n)=(n-b_n)\log(1-g(n,b_n))\leq0$ and $b_n\leq n$,
    \[H(n)\geq n\log(1-g(n,b_n))=-n\sum_{k=1}^\infty\frac{g(n,b_n)^k}{k}.\]
    We already showed in \cref{eqn:n-delta} that $g(n,b_n)\leq n^{-(2+\gamma)}$. 
    Hence
    \[H(n)\geq-\sum_{k=1}^\infty\frac{n^{-k(2+\gamma)+1}}{k}
    \geq -\sum_{k=1}^\infty \frac{n^{-k(1+\gamma)}}{k}=\log(1-n^{-(1+\gamma)})\]
    and so $n(1-e^{H(n)})\leq n^{-\gamma}$. It follows that
     \[np\left[1-\left(1-\Big(1-\frac{b_n}{n-1}\Big)^{b_n}\right)^{n-b_n}\right]=np(1-e^{H(n)})\to0\]
     as $n\to\infty$. We have shown that each term of $|n-\E_{b_n}[X_1^D]|$ converges to $0$ and so $|n-\E_{b_n}[X_1^D]|\to 0$.
\end{proof}

Considering the SARPZF chain $(X_t^S)$, notice that 
\[(1-p)n-\E_b[X_1^S]=(1-p)(n-b)\left(1-\frac{b}{n-1}\right)^b.\]
Notice also that $\E_b[X_1^S]\leq \E_n[X_1^S]=(1-p)n$ by \cref{OneStepExact}.
Following the proof of \cref{OneStepThreshold} through \cref{eqn:n-delta}, one gets the following corollary which gives a threshold for the expected number of blue vertices in SARPZF to be close to $(1-p)n$.
\begin{corollary}
    Let $(X_t^S)$ be the SARPZF Markov chain on $K_n$ with $n\geq 2$ and reversion probability $p\in(0,1)$. If $b_n\geq \sqrt{n\log n^{1+\gamma}}$ with $\gamma>0$, then 
    $$\lim_{n\to\infty}|(1-p)n-\E_{b_n}[X_1^S]|=0.$$
\end{corollary}
In other words, if $K_n$ has $b_n=\Omega(\sqrt{n\log n})$ blue vertices, then with constant $C_1>\sqrt{2}$ we expect $n$ blue vertices after one application of the DARPZF color change rule, and with constant $C_2>1$ we expect $(1-p)n$ blue vertices after one application of the SARPZF color change rule.
It turns out $\sqrt{n\log n}$ is the threshold for this behavior in RPZF. 
Indeed, if $b_n=O(\sqrt{n\log n})$ with constant $C<1$, then the expected number of blue vertices after one step does not converge as in \cref{OneStepThreshold} or its corollary. 
Instead, the SARPZF and DARPZF Markov chains on $K_n$ converge to each other while getting arbitrarily far from $n$.
This is made precise in the next result.
\begin{proposition}\label{SAR-DAR-equal}
    If $b_n\leq \sqrt{n\log n^{1-\gamma}}$ then for any $\varep>0$ there exists an $N$ such that for any $n>N$,
    $|\E_{b_n}[X_1^S]-\E_{b_n}[X_1^D]|<\varep.$ 
    Moreover, $|(1-p)n-\E_{b_n}[X_1^S]|\to\infty$ and $|n-\E_{b_n}[X_1^D]|\to\infty$ as $n\to\infty$.
\end{proposition}
\begin{proof}
    Let $b_n\leq\sqrt{n\log n^{1-\gamma}}$ with $\gamma>0$ and let $g(n,b_n)=\left(1-\frac{b_n}{n-1}\right)^{b_n}$. 
    By \cref{OneStepExact}, to prove $|\E_{b_n}[X_1^S]-\E_{b_n}[X_1^D]|<\varep$ for sufficiently large $n$, it suffices to show that $np(1-g(n,b_n))^{n-b_n}$ converges to $0$ as $n\to\infty$. 
    Observe
    \[g(n,b_n)=e^{b_n\log\left(1-\frac{b_n}{n-1}\right)}
    =e^{b_n\left(-\frac{b_n}{n-1}+O\left(\frac{b_n^2}{n^2}\right)\right)}
    = e^{b_n\left(-\frac{b_n}{n}-\frac{b_n}{n^2-n}+O\left(\frac{b_n^2}{n^2}\right)\right)}
    =e^{-\frac{b_n^2}{n}}e^{O\left(\frac{b_n^3}{n^2}\right)}\]
    since $-\frac{b_n}{n^2-n}=O(b_n^2/n^2)$. 
    Using the expansion $e^x=\sum_{i=0}^\infty x^i/i!$ this is equal to
    \[e^{-\frac{b_n^2}{n}}\left[1+O\left(\tfrac{b_n^3}{n^2}\right)+O\left(\tfrac{b_n^6}{n^4}\right)+\cdots\right]
    =e^{-\frac{b_n^2}{n}}\left[1+O\left(\tfrac{b_n^3}{n^2}\right)\right].\]
    Now apply the assumption $b_n\leq \sqrt{n\log n^{1-\gamma}}$ to get
    $g(n,b)\geq \frac{1}{n^{1-\gamma}}\left[1+O\left(\tfrac{b_n^3}{n^2}\right)\right].$

    Turning to $(1-g(n,b_n))^{n-b_n}$, this can be written as
    \[\exp{-(n-b_n)\sum_{k\geq 1}\frac{g(n,b_n)^k}{k}}\leq e^{-(n-b_n)g(n,b_n)},\]
    and substituting $-g(n,b)\leq -\frac{1}{n^{1-\gamma}}\left[1+O\left(\tfrac{b_n^3}{n^2}\right)\right]$ gives 
    $$(1-g(n,b_n))^{n-b_n}\leq \exp{-(n-b_n)n^{-(1-\gamma)}\left[1+O\left(\tfrac{b_n^3}{n^2}\right)\right]}.$$ 
    Finally, apply $b_n\leq \sqrt{n\log n^{1-\gamma}}$ and simplify to find
    \[0\leq (1-g(n,b_n))^{n-b_n}\leq \exp{-n^\gamma\left(1-\frac{\sqrt{\log n^{1-\gamma}}}{\sqrt{n}}\right)\left[1+O\left(\tfrac{b_n^3}{n^2}\right)\right]}.\]
    Since $\sqrt{\log n^{1-\gamma}/n}\to0$ and $O(b_n^3/n^2)\to0$ as $n\to\infty$, we conclude that $(1-g(n,b_n))^{n-b_n}=O(e^{-n^\gamma})$ from which it follows that $np(1-g(n,b_n))^{n-b_n}\to0$ as $n\to\infty$.
    
    Additionally, $$|n-\E_{b_n}[X_1^D]|\geq(1-p)(n-b_n)g(n,b_n)=|(1-p)n-\E_{b_n}[X_1^S]|$$ and from $g(n,b)\geq \frac{1}{n^{1-\gamma}}\left[1+O\left(\tfrac{b_n^3}{n^2}\right)\right]$ it follows that
    \begin{align*}
        (n-b_n)g(n,b_n)\geq\left(n^\gamma-\frac{b_n}{n^{1-\gamma}}\right)\left[1+O\left(\tfrac{b_n^3}{n^2}\right)\right]
        &\geq \left(n^\gamma-\frac{\sqrt{n\log n^{1-\gamma}}}{n^{1-\gamma}}\right)\left[1+O\left(\tfrac{b_n^3}{n^2}\right)\right]\\
        &=n^\gamma\left(1-\frac{\sqrt{n\log n^{1-\gamma}}}{n}\right)\left[1+O\left(\tfrac{b_n^3}{n^2}\right)\right]
    \end{align*}
    which tends to infinity as $n\to\infty$.
\end{proof}

This result, combined with \cref{OneStepThreshold}, shows that $\sqrt{n\log n}$ is the threshold number of blue vertices for DARPZF on the complete graph to fully force in one step. 
\begin{theorem}\label{thm:threshold}
    Let $(X_t^D)$ be the DARPZF Markov chain on $K_n$ with $n\geq 2$ and reversion probability $p\in(0,1)$. 
    \begin{itemize}
        \item If $b_n\leq \sqrt{n\log n^{1-\gamma}}$ then \(|n-\E_{b_n}[X_1^D]|\to\infty\) as $n\to\infty$, and
        \item if $b_n\geq\sqrt{n\log n^{2+\gamma}}$ then \(|n-\E_{b_n}[X_1^D]|\to0\) as $n\to\infty$.
    \end{itemize}
    That is, $\sqrt{n\log n}$ is a threshold function for expecting DARPZF to fully force in one step.
\end{theorem}
In fact, if $b_n\geq \sqrt{n\log n^{2+\gamma}}$, then since $n\geq X_1^D\geq 0$ we have $|n-\E_{b_n}[X_1^D]|=\E_{b_n}[|n-X_1^D|]\to 0$ as $n\to\infty$. Formally, one says that if $b_n\geq \sqrt{n\log n^{2+\gamma}}$ then $X_1^D$  \textit{converges in mean} to $n$.

The next threshold result, as a consequence of previous, states that if $K_n$ has asymptotically greater than $\sqrt{n\log n}$ blue vertices, then with high probability $K_n$ will be entirely blue after one step of the DARPZF Markov chain. If the number of blue vertices is asymptotically below $\sqrt{n\log n}$, then with high probability the DARPZF chain will not have $n$ blue vertices.
\begin{theorem}\label{KnProbThreshold}
    Let $(X_t^D)$ be the DARPZF Markov chain on $K_n$ with reversion probability $p\in(0,1)$, and let $\gamma>0$.
    \begin{itemize}
        \item If $b_n\leq \sqrt{n\log n^{1-\gamma}}$ then $\P_{b_n}[X_1^D=n]\to0$ as $n\to\infty$, and
        \item if $b_n\geq\sqrt{n\log n^{1+\gamma}}$, then $\P_{b_n}[X_1^D=n]\to1$ as $n\to\infty$.
    \end{itemize}
\end{theorem}
\begin{proof}
    Let $b_n$ be such that $0\leq b_n\leq n$ for all $n$. We may assume $b_n\leq n-2$ since $\P_{n-1}[X_1^D=n]=1$. 
    Notice $\P_{b_n}[X_1^D=n]=q(n,b_n)^{n-b_n}=(1-g(n,b_n))^{n-b_n}$. 
    If $b_n\leq \sqrt{n\log n^{1-\gamma}}$, then in \cref{SAR-DAR-equal} we showed $(1-g(n,b_n))^{n-b_n}=O(e^{-n^\gamma})$. Thus $\P_{b_n}[X_1^D=n]\to 0$ as $n\to\infty$.
    
    Suppose now $b_n\geq\sqrt{n\log n^{1+\gamma}}$. Define $H(n)=(n-b_n)\log\left(1-g(n,b_n)\right)$ so that 
    \[1-\big(1-g(n,b_n)\big)^{n-b_n}=1-e^{H(n)}.\]
    It suffices to show that $1-e^{H(n)}\to 0$ as $n\to\infty$. Equivalently, we show $H(n)\geq \log(1-n^{-\gamma})$ for sufficiently large $n$. The proof of this is almost identical to that in \cref{OneStepThreshold}. Indeed, notice that because $H(n)=(n-b_n)\log(1-g(n,b_n))<0$,
    \[H(n)>n\log(1-g(n,b_n))=-n\sum_{k=1}^\infty\frac{g(n,b_n)^k}{k}.\]
    We showed in \cref{eqn:n-delta} that if $b_n\geq\sqrt{n\log n^{1+\gamma}}$ then $g(n,b_n)\leq n^{-(1+\gamma)}$. 
    Hence
    \[H(n)>-\sum_{k=1}^\infty\frac{n^{-k(1+\gamma)+1}}{k}
    \geq -\sum_{k=1}^\infty \frac{n^{-k\gamma}}{k}=\log(1-n^{-\gamma})\]
    and so $1-(1-g(n,b_n))^{n-b_n}\to0$, which implies $\P_{b_n}[X_1^D=n]=(1-g(n,b_n))^{n-b_n}\to 1$ as $n\to\infty$.
\end{proof}

Since $n\geq X_1^D$, this is equivalent to saying that $X_1^D$ \textit{converges in probability} to $n$ when $b_n\geq \sqrt{n\log n^{1+\gamma}}$. Formally,
\(\P_{b_n}[|n-X_1^D|>\varep]\to0\)
for all $\varep>0$ which is equivalent to
\(\P_{b_n}[|n-X_1^D|=0]\to1,\)
and so $X_1^D=n$ with high probability.
Finally, we establish when the upper bound presented in \cref{lem:OneStepUpperbound} is asymptotically tight.

\begin{proposition}
    Let $(X_t)$ be an RPZF Markov chain on $K_n$ with reversion probability $p$. If $b_n\leq\frac{\sqrt{n}}{\log n}$, then for any $\varep>0$, 
    $$(1-p)(b_n+(1-\varep){b_n}^2)\leq \E_{b_n}[X_1]\leq (1-p)(b_n+{b_n}^2)$$
    for $n$ sufficiently large. In particular, if $b\in\naturals$ is fixed then $\E_b[X_1]\to(1-p)(b+b^2)$ as $n\to\infty$. 
\end{proposition}
\begin{proof}
    Since $b_n\leq\frac{\sqrt{n}}{\log n}\leq\sqrt{n\log n}$, by \cref{SAR-DAR-equal} we need consider only the SARPZF chain.
    By \cref{lem:OneStepUpperbound}, $\E_{b_n}[X_1]\leq (1-p)(b_n+b_n^2)$.
    Let $F_1$ denote the number of vertices forced at time $t=1$ during phase 1. Then by \cref{OneStepExpected}, $\E_{b_n}[X_1]=(1-p)(b_n+\E_{b_n}[F_1]).$ The authors of \cite{CC+20} showed in the proof of Theorem 3.1 that $\E_{b_n}[F_1]=b_n^2-o(b_n^2)$ when $b_n\leq\frac{\sqrt{n}}{\log n}$. 
    Consequently, since $b_n^2=o(n)$, for any $\varep>0$, $\E_{b_n}[F_1]>(1-\varep)b_n^2$ for sufficiently large $n$. Thus $\E_{b_n}[X_1]>(1-p)(b_n+(1-\varep)b_n^2)$ and for fixed $b\in\naturals$, taking $\varep\to0$ gives that $\E_b[X_1]\to (1-p)(b+b^2)$ as $n\to\infty$.
\end{proof}

\subsection{The balanced complete bipartite graph}
The \defn{complete bipartite graph} $K_{m,n}$ is the graph of order $m+n$ whose vertices can be partitioned into two parts $U=\{u_1,\ldots,u_m\}$ and $V=\{v_1,\ldots,v_n\}$ such that the edges of the graph are $u_iv_j$ for all $1\leq i\leq m$ and all $1\leq j\leq n$. 
If $m=n$ then $K_{n,n}$ is the \defn{balanced} complete bipartite graph. 
We will see in \cref{sec:sims} that $K_{n,n}$ behaves very similarly to $K_{2n}$ in DARPZF. 
This section develops some partial characterizations of the behavior of $K_{n,n}$ which support this observation.

\begin{lemma}\label{KnnProb}
    Let $K_{m,n}$ have vertex partitions $U$ and $V$. Suppose $b_U$ vertices in $U$ are blue and $b_V$ vertices in $V$ are blue. The the probability that $U$ forces $V$ entirely blue in one step is
    \[\P[U\to V]=\left(1-\left(1-\frac{b_V+1}{|V|}\right)^{b_U}\right)^{|V|-b_V}.\]
\end{lemma}

\emph{Proof.}
    Let $B$ be the set of blue vertices in $K_{m,n}$. Let $u\in U\cap B$ be blue and $v\in V\setminus (B\cap V)$ be white. Then 
    \[\P[u\to v]=\frac{|N[u]\cap B|}{\deg u}=\frac{b_V+1}{|V|}\]
    and so the probability that $v$ is forced by some vertex in $U$ is
    \begin{align*}
        \P[U\to v]=1-\P[U\not\to v]
        =1-\prod_{u\in U\cap B}\P[u\not\to v]
        =1-\left(1-\frac{b_V+1}{|V|}\right)^{b_U}.
    \end{align*}
    Thus
    \[\P[U\to V]=\prod_{v\in V\setminus(B\cap V)}\P[U\to v]=\left(1-\left(1-\frac{b_V+1}{|V|}\right)^{b_U}\right)^{|V|-b_V}.~~~~~\cvd\]

We now give an upper bound for the threshold number of blue vertices to fully force the balanced complete bipartite in one step with high probability, starting from one of two cases.
The first case is when the vertex parts $U$ and $V$ of $K_{n,n}$ have the same number of blue vertices, and the second is when, without loss of generality, $U$ is entirely blue and $V$ is minimally blue. 
\begin{proposition}\label{Knn1}
    Let $K_{n,n}$ have vertex parts $U$ and $V$ with at least $b_n^U$ and $b_n^V$ blue vertices in $U$ and $V$, respectively. 
    \begin{itemize}
        \item If $b_n^U\geq \sqrt{n\log n^{1+\gamma_1}}$ and $b_n^V\geq \sqrt{n\log n^{1+\gamma_2}}$ for any $\gamma_1,\gamma_2>0$, then with high probability $K_{n,n}$ is blue after one application of the DARPZF color change rule.
        \item If $b_n^U=n$ and $b_n^V\geq\log(n^{1+\gamma})$ for any $\gamma>0$, then with high probability $K_{n,n}$ is blue after one application of the DARPZF color change rule.
    \end{itemize}
\end{proposition}
\begin{proof}
    Suppose first $b_n^U\geq \sqrt{n\log n^{1+\gamma_1}}$ and $b_n^V\geq \sqrt{n\log n^{1+\gamma_2}}$ for some $\gamma_1,\gamma_2>0$.
    Let $\gamma=\min\{\gamma_1,\gamma_2\}$.
    We may assume $U$ and $V$ each have $b_n\geq \sqrt{n\log n^{1+\gamma}}$ blue vertices because the probability of fully forcing in one step monotonically increases in both $b_n^V$ and $b_n^U$.
    Since the events $\{U\to V\}$ and $\{V\to U\}$ are independent (when both occuring at time $t$), the probability that $K_{n,n}$ is blue after one step of DARPZF is
    \[\P[U\to V]\P[V\to U]=\left(\left(1-\left(1-\frac{b_n+1}{n}\right)^{b_n}\right)^{n-b_n}\right)^2=\left(1-\left(1-\frac{b_n+1}{n}\right)^{b_n}\right)^{2n-2b_n}\]
    by \cref{KnnProb}.
    Call this probability $P(b_n)$. We wish to show $P(b_n)\to 1$ as $n\to\infty$. 
    Note that $P(b_n)=1$ if $b_n=n$ so assume $b_n\leq n-1$.
    Define 
    $f(n,b_n)=\left(1-\frac{b_n+1}{n}\right)^{b_n}$
    and let $H(n)=(2n-2b_n)\log(1-f(n,b_n))$. Then $P(b_n)=e^{H(n)}$ and so $P(b_n)\to 1$ if $H(n)\to 0$. Observe that $H(n)<0$ and hence
    \[H(n)>2n\log(1-f(n,b_n))=-2n\sum_{k= 1}^\infty\frac{f(n,b_n)^k}{k}.\]
    Now in the style of \cref{KnProbThreshold}, if $f(n,b_n)\leq n^{-(1+\gamma)}$ then
    \[H(n)>-2\sum_{k=1}^\infty\frac{n^{-k(1+\gamma)+1}}{k}\geq-2\sum_{k=1}^\infty\frac{n^{-k\gamma}}{k}=2\log(1-n^{-\gamma})\to0\]
    as $n\to\infty$. To see that $f(n,b_n)\leq n^{-(1+\gamma)}$, observe
    \begin{align*}
    f(n,b_n)=\exp{b_n\log\left(1-\frac{b_n+1}{n}\right)}
    &=\exp{-b_n\sum_{k=1}^\infty\frac{(b_n+1)^k}{kn^k}}\\
    &=\exp{-b_n\frac{b_n+1}{n}}\exp{-b_n\sum_{k=2}^\infty\frac{(b_n+1)^k}{kn^k}}
    \end{align*}
    which is bounded above by $e^{-b_n^2/n}$ since $b_n>0$. Then because $b_n\geq\sqrt{n\log n^{1+\gamma}}$, 
    \[f(n,b_n)\leq e^{-b_n^2/n}\leq n^{-(1+\gamma)}.\]
    Hence $H(n)<0$ implies $H(n)\to0$ and thus $P(b_n)=e^{H(n)}\to 1$ as $n\to\infty$.

    Now assume $b_n^U=n$ and $b_n^V\geq\log(n^{1+\gamma})$ for some $\gamma>0$. 
    Then the probability that $K_{n,n}$ is entirely blue after one step of DARPZF is
    \[\P[U\to V]=\left(1-\left(1-\frac{b_n^V+1}{n}\right)^n\right)^{n-b_n^V}.\]
    Let $f(n,b_n^V)=\left(1-\frac{b_n^V+1}{n}\right)^{n}$ and $H(n)=(n-b_n^V)\log(1-f(n,b_n^V))$. Notice that $f$ differs from before in the exponential. Like before, it suffices to show that $f(n,b_n^V)\leq n^{-(1+\gamma)}$ because then
    \[H(n)>n\log(1-f(n,b_n^V))=-n\sum_{k=1}^\infty\frac{f(n,b_n^V)}{k}\geq -\sum_{k=1}^\infty\frac{n^{-k\gamma}}{k}=\log(1-n^{-\gamma}).\]
    Now, 
    \begin{align*}
        f(n,b_n^V)=\exp{n\log\left(1-\frac{b_n^V+1}{n}\right)}
        &=\exp{-n\sum_{k=1}^\infty\frac{(b_n^V+1)^k}{kn^k}}\\
        &=\exp{-n\frac{b_n^V+1}{n}}\exp{-n\sum_{k=2}^\infty\frac{(b_n^V+1)^k}{kn^k}}
        \end{align*}
    which is bounded above by $e^{-b_n^V}$ since $b_n^V>0$. Then because $b_n^V\geq\log n^{1+\gamma}$ it follows that \(f(n,b_n^V)\leq n^{-(1+\gamma)}.\) Thus, $0>H(n)>\log(1-n^{-\gamma})\to0$ and so $\P[U\to V]=e^{H(n)}\to 1$ as $n\to\infty$.
\end{proof}
This result supports the claim that when $n$ is large and the number of blue vertices is balanced, $K_{n,n}$ behaves similarly to $K_{2n}$.

\subsection{The star graph}
The star graph on $n$ vertices is $K_{1,n-1}$ and the singleton vertex is called the \defn{universal vertex} because it is adjacent to all other vertices.
In this section we show that the star graph exhibits a large threshold value for one-step forcing.
Let $K_{1,n-1}$ be the star graph on $n$ vertices with universal vertex $v$ and set of currently blue vertices $B$. 
Notice that if $|B|=b\leq n-2$, the only way for $K_{1,n-1}$ to be forced in one step is if $v\in B$ because $v$ is the only neighbor of white degree 1 vertices. 
Hence, when calculating the one-step threshold for $K_{1,n-1}$ we need consider only when $v\in B$. 
In that case, $K_{1,n-1}$ is fully forced in one step with probability 
\[\P[v\rightarrow V(K_{1,n-1})\setminus B]=\left(\frac{b}{n-1}\right)^{n-b}.\]
This lets us calculate the one-step expected values for the star.
\begin{lemma}
Let $(X_t^S)$ and $(X_t^D)$ be the SARPZF and DARPZF Markov chains on $K_{1,n-1}$ with universal vertex $v$. If $v$ blue at time $t=0$, then
\[\E_b[X_1^S]=(1-p)\left(b+(n-b)\left(\frac{b}{n-1}\right)\right)\]
and
\[\E_b[X_1^D]=\E_b[X_1^S]+np\left(\frac{b}{n-1}\right)^{n-b}.\]
\end{lemma}

{\em Proof.}
Let $w_1,\ldots,w_{n-b}$ denote the white vertices. Then by \cref{OneStepExpected}, 
\[\E_b[X_1^S]=(1-p)\big(b+(n-b)\P_b[v\to w_i]\big)=(1-p)\left(b+(n-b)\left(\frac{b}{n-1}\right)\right).\]
Using the same approach as that in \cref{OneStepExact}, one also calculates
\[\E_b[X_1^D]=\E_b[X_1^S]+np\left(\frac{b}{n-1}\right)^{n-b}.~~~~~~\cvd\]

Observe that if $K_{1,n-1}$ has $b_n=n-1-C$ vertices blue, then 
$$\P[v\rightarrow V(K_{1,n-1})\setminus B]=\left(1-\frac{C}{n-1}\right)^{1+C}\to1$$
as $n$ grows to infinity. This turns out to be the threshold for fully forcing in DARPZF: if the distance between $n$ and $b_n$ is unbounded, then with high probability $K_{1,n-1}$ is not fully forced in the next step.
\begin{theorem}\label{starThm}
    Let $(X_t^D)$ be the DARPZF Markov chain on $K_{1,n-1}$ with universal vertex $v$ and reversion probability $p\in(0,1)$. When $v$ is blue, we have the following:
    \begin{itemize}
        \item if $b_n=n-1-C$ for some constant $C\in\naturals$, then \(|n-\E_{b_n}[X_1^D]|\to pC(C+1)\) as $n\to\infty$, and
        \item if $b_n= n-1-\omega(1)$, then \(|n-\E_{b_n}[X_1^D]|\to\infty\) as $n\to\infty$.
    \end{itemize}
\end{theorem}
\begin{proof}
Let $(X_t^D)$ be the DARPZF Markov chain on $K_{1,n-1}$ with reversion probability $p\in(0,1)$, let $0\leq b_n\leq n-1$, and let $v$ denote the universal vertex. 
Assume $v$ is blue.
If we replace $g(n,b_n)=\left(1-\frac{b_n}{n-1}\right)^{b_n}$ with $1-\frac{b_n}{n-1}$ in the proof of \cref{OneStepThreshold}, then we may simplify $n-\E_{b_n}[X_1^D]$ 
as
\begin{align}\label{eqn:starexp}
n-\E_{b_n}[X_1^D]
&=(1-p)(n-b_n)\left(1-\frac{b_n}{n-1}\right)+np\left[1-\left(\frac{b_n}{n-1}\right)^{n-b_n}\right].
\end{align}

Suppose $b_n=n-1-C$ for some constant $C\geq 0$. Then
\begin{align*}
    n-\E_{b_n}[X_1^D]&=(1-p)(n-(n-1-C))\left(1-\frac{n-1-C}{n-1}\right)+np\left(1-\left(\frac{n-1-C}{n-1}\right)^{n-(n-1-c)}\right)\\
    &=(1-p)(1+C)\frac{C}{n-1}+p\left(n-n\left(1-\frac{C}{n-1}\right)^{C+1}\right).
\end{align*}
It is immediate that $(1-p)(1+C)\frac{C}{n-1}\to0$ as $n\to\infty$. To see that $p\left(n-n\left(1-\frac{C}{n-1}\right)^{C+1}\right)\to pC(C+1)$, observe
\begin{align*}
    n-n\left(1-\frac{C}{n-1}\right)^{C+1}&=n-n\sum_{k=0}^{C+1}(-1)^k\binom{C+1}{k}\left(\frac{C}{n-1}\right)^k\\
    &=n-n+n(C+1)\frac{C}{n-1}-n\sum_{k=2}^{C+1}(-1)^{k}\binom{C+1}{k}\left(\frac{C}{n-1}\right)^k\\
    &=C(C+1)\frac{n}{n-1}+O(n^{-1}).
\end{align*}
Taking $n\to\infty$ gives $n-\E_{b_n}[X_1^D]\to pC(C+1)$.

On the other hand, suppose $b_n=n-1-f_n$ with $f_n=\omega(1)$. We want to show 
\[n-\E_{b_n}[X_1^D]=(1-p)(n-b_n)\left(1-\frac{b_n}{n-1}\right)+np\left[1-\left(\frac{b_n}{n-1}\right)^{n-b_n}\right]\to\infty\]
If $f_n=\omega(\sqrt{n})$ then 
\[(1-p)(n-b_n)\left(1-\frac{b_n}{n-1}\right)=(1-p)(1+f_n)\frac{f_n}{n-1}=(1-p)\frac{f_n+f_n^2}{n-1}\to\infty,\]
so assume $f_n=O(\sqrt{n})$.
Notice
\[np\left[1-\left(\frac{b_n}{n-1}\right)^{n-b_n}\right]=np\left[1-\left(1-\frac{f_n}{n-1}\right)^{1+f_n}\right]\]
and so define $h_n(x)=x^{1+f_n}$.
Now $h_n'(x)=(1+f_n)x^{f_n}$ and then applying Taylor's theorem to $h_n(x)$ around $x_0=1$,
\[h_n(x)=h_n(1)+h_n'(\xi)(x-1)\]
for some $\xi\in(x,1)$. Thus
\[h_n\left(1-\frac{f_n}{n-1}\right)=1-(1+f_n)\xi_n^{f_n}\frac{f_n}{n-1}\]
for some $\xi_n\in\left(1-\frac{f_n}{n-1},1\right)$ and it follows that 
\begin{align*}
np\left[1-\left(\frac{b_n}{n-1}\right)^{n-b_n}\right]
=np\left[1-h_n\left(1-\frac{f_n}{n-1}\right)\right]
=np\left[(1+f_n)\xi_n^{f_n}\frac{f_n}{n-1}\right].
\end{align*}
which simplifies to $\frac{np}{n-1}\xi_n^{f_n}(f_n+f_n^2).$
By definition, $f_n+f_n^2\to\infty$ when $f_n=\omega(1)$. 
However, $\xi_n^{f_n}$ must be accounted for; if we can show $\xi_n^{f_n}\geq M$ for some constant $M>0$ not dependent on $n$, then $\frac{np}{n-1}\xi_n^{f_n}(f_n+f_n^2)\geq \frac{np}{n-1}M(f_n+f_n^2)\to\infty$ as desired.
Recall $\xi_n\in\left(1-\frac{f_n}{n-1},1\right)$ and so $\xi_n^{f_n}>\left(1-\frac{f_n}{n-1}\right)^{f_n}$. 
Observe
\[\left(1-\frac{f_n}{n-1}\right)^{f_n}=\exp{f_n\log\left(1-\frac{f_n}{n-1}\right)}\] and using the Taylor expansion $\log(1-x)=-\sum_{k=1}^\infty x^k/x$ for $|x|<1$,
\begin{align}\label{eqn:fn}
f_n\log\left(1-\frac{f_n}{n-1}\right)=-f_n\sum_{k\geq 1}\frac{1}{k}\left(\frac{f_n}{n-1}\right)^k\nonumber
=-\sum_{k\geq 1}\frac{f_n^{k+1}}{k(n-1)^k}.
\end{align}
It follows from $f_n=O(\sqrt{n})$ that
\[f_n\log\left(1-\frac{f_n}{n-1}\right)=O(1)\]
and so there exists an $M>0$ independent of $n$ such that for sufficiently large $n$, 
$\left|f_n\log\left(1-\frac{f_n}{n-1}\right)\right|<M.$
Equivalently, $\left(1-\frac{f_n}{n-1}\right)^{f_n}>e^{-M}$. 
Thus $\xi_n\in\left(1-\frac{f_n}{n-1},1\right)$ implies $\xi_n^{f_n}>\left(1-\frac{f_n}{n-1}\right)^{f_n}> e^{-M}>0$ and so
\[np\left[1-\left(\frac{b_n}{n-1}\right)^{n-b_n}\right]=\frac{2np}{n-1}\xi_n^{f_n}(f_n+f_n^2)\to\infty\]
since $f_n=\omega(1)$.

It is left to consider when $f_n\neq O(\sqrt{n})$ and $f_n\neq\omega(\sqrt{n})$. In this case, define $g_n=\sup_{m\geq n} f_m$. Then $g_n\geq f_n$ and thus $n-1-f_n\geq n-1-g_n$. Since $\E_{b_n}[X_1^D]$ monotonically increases in $b_n$ (keeping $v$ blue), it suffices to show that $n-\E_{n-1-g_n}[X_1^D]\to\infty$. We claim that $g_n=\omega(\sqrt{n})$. Indeed, let $M>0$ and $N\in\naturals$ be arbitrary. Because $f_n\neq O(\sqrt{n})$, there exists an $n_0>N$ such that $f_{n_0}>M\sqrt{n_0}$. Then
\[g_N=\sup_{m\geq N} f_m\geq f_{n_0}>M\sqrt{n_0}>M\sqrt{N}\]
and so $g_n=\omega(\sqrt{n})$. It follows from before that $n-\E_{n-1-g_n}[X_1^D]\to\infty$ and so $n-\E_{n-1-f_n}[X_1^D]\to\infty$. 
\end{proof}

Observe that $|n-\E_{b_n}[X_1^D]|=\E_{b_n}[|n-X_1^D|]$ since $n\geq\E_{b_n}[X_1^D]$ and so \cref{starThm} says $X_1^D$ converges in mean to $n$ exactly when $b_n=n-1$ for sufficiently large $n$. As a consequence, $X_1^D$ converges to $n$ in probability and hence $X_1^D=n$ with high probability.

\section{DARPZF Simulations and Approximations}\label{sec:sims}
Consider the DARPZF chain on $K_n$ with reversion probability $p\in(0,1)$. 
Using \cref{thm:DieVsAbs} we calculate $p_D(K_n,S_1)$, the reversion probability such that $K_n$ has equal probability to either die out or fully force when starting from a single blue vertex,  for small $n$.

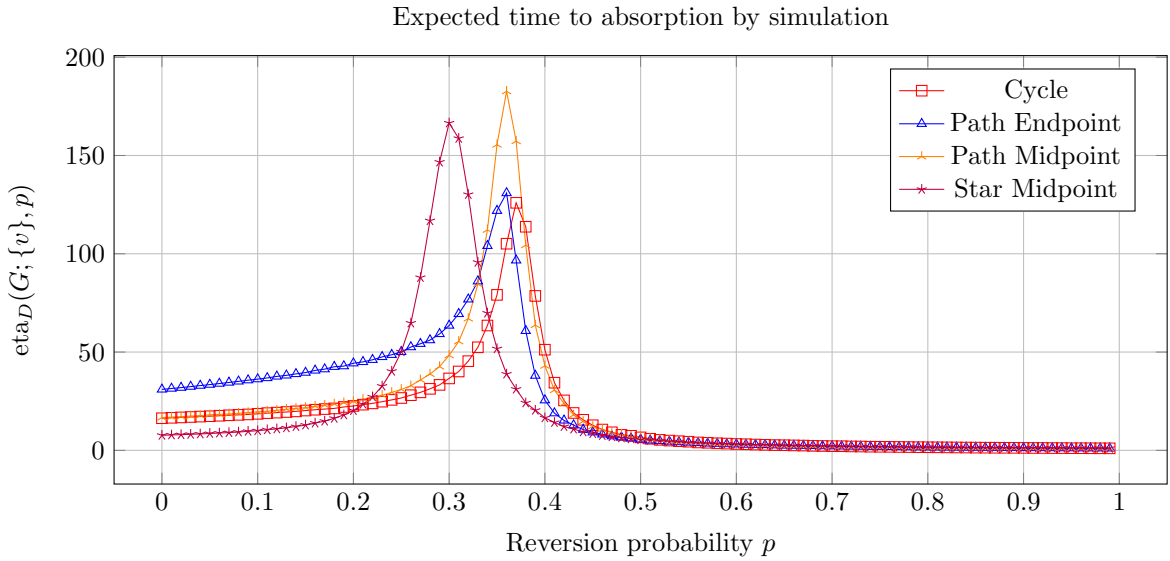
\begin{figure}[h]
\centering
\begin{tikzpicture}
    \begin{axis}[
        width=14cm,
        height=5.7cm,
        xlabel={Reversion probability $p$},
        ylabel={$\expabs_D(G;\{v\},p)$},
        xmin=-0.05,
        xmax=1.05,
        title={Expected time to absorption by simulation},
        grid=major,
        legend pos=north east,
        scale only axis
    ]
        \addplot [color=red, mark=square] table [col sep=comma] {abs_cycle32_10000.txt};
        \addlegendentry{Cycle}

        \addplot [color=blue, mark=triangle] table [col sep=comma] {abs_path32_10000_end.txt};
        \addlegendentry{Path Endpoint}

        \addplot [color=orange, mark=Mercedes star] table [col sep=comma] 
        {abs_path32_10000_mid.txt};
        \addlegendentry{Path Midpoint}

        \addplot [color=purple, mark=star] table [col sep=comma] 
        {abs_star32_10000_mid.txt};
        \addlegendentry{Star Midpoint}
    \end{axis}
\end{tikzpicture}
\caption{Monte Carlo simulations of DARPZF on the cycle, path, and star on 32 vertices, starting from a single vertex.}\label{fig:abs_time}
\end{figure}

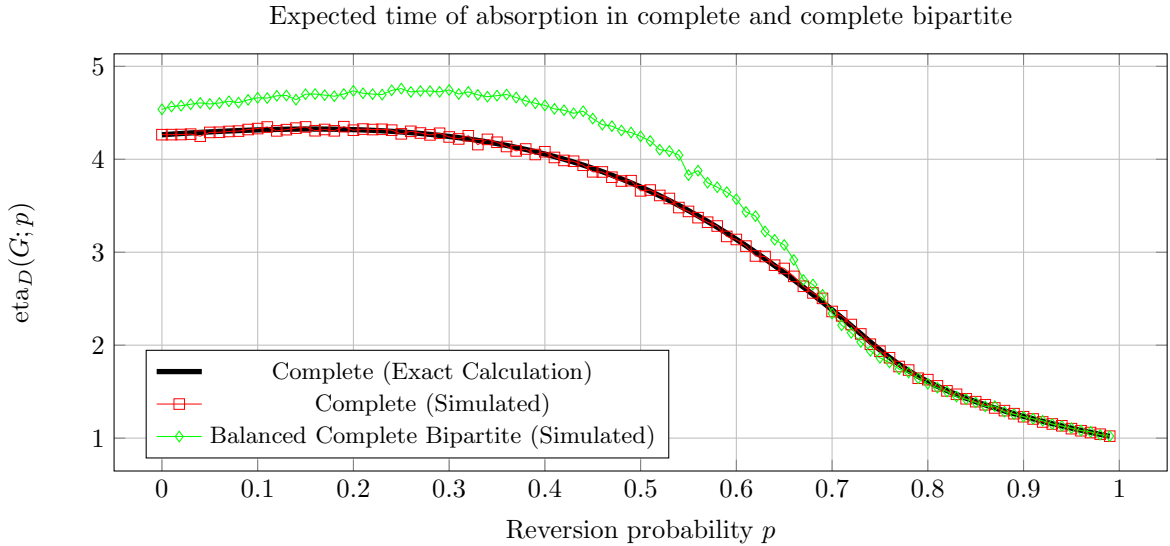
\begin{figure}[h]
\centering
\begin{tikzpicture}
    \begin{axis}[
        legend style={font=\small},
        width=14cm,
        height=5.55cm,
        xmin=-0.05,
        xmax=1.05,
        xlabel={Reversion probability $p$},
        ylabel={$\expabs_D(G;p)$},
        title={Expected time of absorption in complete and complete bipartite},
        grid=major,
        legend pos=south west,
        scale only axis
    ]
        \addplot [color=black, line width=2pt] table [col sep=comma] {abs_K32_exact.txt};
        \addlegendentry{Complete (Exact Calculation)}

        \addplot [color=red, mark=square, line width=0pt] table [col sep=comma] {abs_K32_sim_10000_2.txt};
        \addlegendentry{Complete (Simulated)}

        \addplot [color=green, mark=diamond] table [col sep=comma] {abs_completebi32_10000_end.txt};
        \addlegendentry{Balanced Complete Bipartite (Simulated)}
    \end{axis}
\end{tikzpicture}
\caption{Monte Carlo simulations of DARPZF on $K_{32}$ and $K_{16,16}$ compared to exact calculations for $K_{32}$, starting from a single vertex.}\label{fig:abs_time_Kn}
\end{figure}

\begin{table}[h]
\begin{center}
\begin{tabular}{|c|c||c|c||c|c||c|c|}
\hline
    $n$ & $p_D(K_n,S_1)$ & $n$ & $p_D(K_n,S_1)$ & $n$ & $p_D(K_n,S_1)$& $n$ & $p_D(K_n,S_1)$ \\
    \hline
     $3$ & $.6$ & $8$ & $0.427761$ & $13$&$0.433535$&$18$&$0.435628$\\
     \hline
     $4$ & $0.466548$ & $9$ & $0.429115$ & $14$&$0.434157$& $19$ & $0.43585$ \\
     \hline
     $5$ & $0.437779$& $10$ & $0.43052$ &$15$&$0.434648$ & $20$ & $0.43604$\\
     \hline
     $6$ & $0.428853$ & $11$ & $0.431747$ &$16$&$0.435042$ & $21$ & $0.436203$\\
     \hline
     $7$ & $0.427101$ &$12$ & $0.432745$&$17$&$0.435363$&$22$&$0.436346$\\
     \hline
\end{tabular}
\end{center}
\vspace*{1mm}
\caption{Exact calculations of the critical reversion probability for the complete graph}
\end{table}
\vspace*{-5mm}
These values are calculated from Mathematica's symbolic matrix operations. We can further estimate $p_D(K_n,S_1)$ for larger $n$ to arbitrary levels of precision, as provided in \cref{fig:KnCritical}.
\begin{table}[!h]
\begin{center}
\begin{tabular}{|c|c||c|c|}
\hline
    $n$ & $p_D(K_n,S_1)$ & $n$ & $p_D(K_n,S_1)$ \\
    \hline
     $12$ & $0.43274\pm 0.00001 $ & $96$ & $0.43805\pm0.00005$ \\
     \hline
     $16$ & $0.43505\pm 0.00005$ & $128$ & $0.43815\pm 0.00005$ \\
     \hline
     $32$ & $0.43715\pm 0.00005$& $156$ & $0.43818\pm 0.00005$ \\
     \hline
     $64$ & $0.4379\pm 0.00005$ & $192$ & $0.4382\pm0.00005$ \\
     \hline
\end{tabular}
\end{center}
\caption{Numerically approximated critical reversion probability for the complete graph.}\label{fig:KnCritical}
\end{table}
Going beyond the critical reversion probability for the complete graph, \cref{fig:abs_time,fig:abs_time_Kn} approximate the expected time of absorption for various graphs using Monte Carlo simulations, starting from a single blue vertex. Let $V=\{v_1,\ldots,v_n\}$ be a set of vertices. The \textit{path graph} $P_n$ is the graph with edges $v_iv_{i+1}$ for $1\leq i\leq n-1$. The endpoints are $v_1$ and $v_n$, and the midpoint is $v_{\lceil n/2\rceil}$. The \textit{cycle graph} is a path with the additional edge $v_nv_1$. When simulating DARPZF on the path, starting from an endpoint and starting from a midpoint has a noticeable affect on the expected time to absorption.

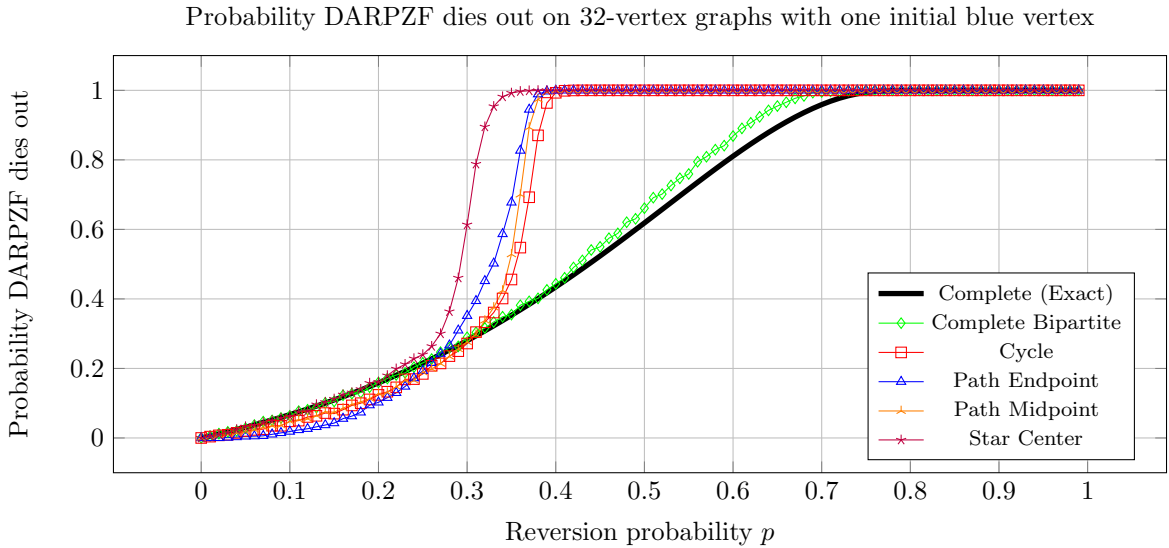
\begin{figure}[!h]
\centering
\begin{tikzpicture}
    \begin{axis}[
        legend style={font=\footnotesize},
        width=14cm,
        height=5.55cm,
        xlabel={Reversion probability $p$},
        ylabel={Probability DARPZF dies out},
        title={Probability DARPZF dies out on 32-vertex graphs with one initial blue vertex},
        grid=major,
        legend pos=south east,
        scale only axis
    ]
        \addplot [color=black, line width=2pt] table [col sep=comma] {Kn32_exact.txt};
        \addlegendentry{Complete (Exact)}
        
        \addplot [color=green, mark=diamond] table [col sep=comma] {completebi32_8000_end.txt};
        \addlegendentry{Complete Bipartite}
        
        \addplot [color=red, mark=square] table [col sep=comma] {cycle32_8000.txt};
        \addlegendentry{Cycle}

        \addplot [color=blue, mark=triangle] table [col sep=comma] {path32_8000_end.txt};
        \addlegendentry{Path Endpoint}

        \addplot [color=orange, mark=Mercedes star] table [col sep=comma] 
        {path32_8000_mid.txt};
        \addlegendentry{Path Midpoint}

        \addplot [color=purple, mark=star] table [col sep=comma] 
        {star32_10000_mid.txt};
        \addlegendentry{Star Center}
    \end{axis}
\end{tikzpicture}
\caption{Monte Carlo simulations of DARPZF on the cycle, path, star, and balanced complete bipartite graphs on 32 vertices, starting from a single blue vertex.}\label{fig:dieout}
\end{figure}

Finally, \cref{fig:dieout} presents Monte Carlo simulations of DARPZF on different 32 vertex graphs starting from one blue vertex. Notice that when the reversion probability $p>0.4$, the cycle, path, and star die out with probability nearly one. On the other hand, the highly connected complete graph and balanced complete bipartite graph don't reach the same level die out until $p>0.75$ and $p>0.7$, respectively.

The critical reversion probability when starting at a single vertex can also be estimated from \cref{fig:dieout}. 
Observe that for the cycle, path, and star graphs, the $p$ which maximizes $\expabs_D(G;p)$ in \cref{fig:abs_time} closely aligns with the approximate $p_D(G,\{v\})$ from \cref{fig:dieout}.
On the other hand, the maximums of $\expabs_D(K_{32};p)$ and $\expabs_D(K_{16,16};p)$ occur at a $p$ noticeably smaller than their respective critical reversion probabilities. 
In the case of the complete graph, we know from \cref{thm:threshold} that only a small percentage of vertices are needed to have a high chance of forcing the entire graph blue. 
Compare this with the cycle or path graph, where the threshold for one-step fully forcing is of order $n$. 
Indeed, at least one third of the vertices must be blue to even have a non-zero probability to fully force the cycle or path in one step.

When trying to explain the difference in RPZF graph behavior, one can also also consider their expected propagation times in traditional PZF. The cycle and path are expected to take $\Theta(n)$ steps to force in PZF \cite{GH18}. Contrast this with the complete graph which has an expected propagation time of $\Theta(\log\log n)$ steps \cite{CC+20}. The cycle and path thus have more time to be stymied by the reversion of blue vertices in RPZF, which may result in longer-lasting oscillations of white and blue vertices. The complete graph, on the other hand, forces significantly more quickly and so does not as easily fall into an oscillation of white and blue vertices.

In RPZF, both the one step threshold and expected propagation time are a means of quantifying how ``hard" a graph is to fully force. 
This is in contrast to traditional PZF where propagation time does not correlate to difficulty. Indeed, consider the path or cycle on $n$ vertices when starting from a single blue vertex.
These graphs require only a single force to become deterministic in nature, but will always take at least $n/2$ rounds to fully force.

When starting from a single vertex, \cref{fig:abs_time_Kn,fig:dieout} indicate that the balanced complete bipartite graph behaves much more like the complete graph than the star graph (a severely ``unbalanced" complete bipartite graph). This observation is supported by \cref{Knn1} which gives the one-step forcing threshold for $K_{n/2,n/2}$ as order $O\left(\sqrt{n\log(n/2)}\right)$, which is comparable to the $K_n$ threshold of $\Theta\left(\sqrt{n\log n}\right)$ and significantly smaller than the $K_{1,n-1}$ threshold of $n-1-o(1)$ from \cref{starThm}.

\appendix

\section{Background Analysis and Probability}\label{sec:appendix}
Let $f(n)$ and $g(n)$ be functions from the nonnegative integers to the real numbers, where $g$ is strictly positive for sufficiently large input.
Write $f=O(g)$ if there exists constants $C,N>0$ such that for all $n>N$, $|f(n)|\leq Cg(n)$ and write $f=o(g)$ if for all $C>0$ there exists an $N>0$ such that for all $n>N$, $|f(n)|<Cg(n)$. Symmetrically, $f=\Omega(g)$ if there exists constants $C,N$ such that for all $n>N$, $f(n)\geq Cg(n)$ and $f=\omega(g)$ if for all $C>0$ there exists an $N>0$ such that for all $n>N$, $f(n)>Cg(n)$. 
That is, $f=O(g)$ if and only if $g=\Omega(f)$ and $f=o(g)$ if and only if $g=\omega(f)$.
One can also define this notions in terms of limit behavior. Namely, $f=o(g)$ when
\(\lim_{n\to\infty}\frac{f(n)}{g(n)}=0\)
and $f=\omega(g)$ if
\(\lim_{n\to\infty}\frac{f(n)}{g(n)}=\infty.\)
We refer the reader to  e.g. \cite{de1981asymptotic} for a more thorough introduction to asymptotic notation.

We shall make use of the following standard facts about asymptotic notation.  Let $f$ and $g$ be real-valued functions.
\begin{fact}\label{AsymptLemma}
    If $f=o(g)$ then $f=O(g)$, and if $f=\omega(g)$ then $f=\Omega(g)$. Moreover, $f=O(f)$ and $f=\Omega(f)$. 
\end{fact}

\begin{fact}
    $O(f)O(g)=O(fg)$ and $\Omega(f)\Omega(g)=\Omega(fg).$
    Moreover, $o(f)O(g)=o(fg)$ and $\omega(f)\Omega(g)=\omega(fg)$. 
\end{fact}

\begin{fact}
    If $f=O(g)$, then $O(f)+O(g)=O(g)$. 
\end{fact}

Fubini's theorem gives a condition for switching the order of integration in a double integral. As a consequence, is it gives a condition for interchanging the order of infinite sums.
\begin{theorem}[Fubini's theorem]\label{fubini}
    Let $f:\reals^2\to R$ be measurable. If 
    \[\int_{\reals^2}|f|<\infty\]
    then
    \[\int_{\reals}\left(\int_{\reals}f(x,y)\,\textnormal{d}x\right)\textnormal{d}y
    =\int_{\reals}\left(\int_{\reals}f(x,y)\,\textnormal{d}y\right)\textnormal{d}x
    =\int_{\reals^2}f.\]
    In particular, if $f:\naturals^2\to\reals$ is such that 
    \[\sum_{(n,m)\in\naturals^2} |f(n,m)|<\infty,\]
    then
    \[\sum_{n=0}^\infty\sum_{m=0}^\infty f(n,m)=\sum_{m=0}^\infty\sum_{n=0}^\infty f(n,m)=\sum_{(n,m)\in\naturals^2} f(n,m).\]
\end{theorem}

\bigskip
{\bf Acknowledgment.} This research was partially supported by NSF grant DMS-1839918 and the author thanks the National Science Foundation for this support.


\end{document}